\documentclass[12pt]{article}

\usepackage[left=1.5in, right=1in, top=1in, bottom=1in, dvips,pdftex]{geometry}
\usepackage{color, setspace, graphicx, verbatim, mathrsfs, amsmath, amsthm, amssymb, amsfonts,  multirow, moreverb}
\newtheorem{defn}{Definition}[section]

\newtheorem{thm}[defn]{Theorem}

\theoremstyle{remark}

\numberwithin{equation}{section} \numberwithin{figure}{subsection}

\DeclareMathOperator{\airy}{Ai}

\def\ra{\rightarrow}
\def\iy{\infty}

\def\be{\begin{equation}}
\def\ee{\end{equation}}

\newcommand{\bE}{\mathbb{E}}
\newcommand{\bP}{\mathbb{P}}

\begin{document}

\title{\textbf{Edgeworth Expansion of the Largest Eigenvalue Distribution Function
 of GOE }}

\author{Leonard N.~Choup \\ Department of Mathematical Sciences \\ University Alabama in Huntsville\\
 Huntsville, AL 35899, USA \\
 email:  \texttt{Leonard.Choup@uah.edu}
} \maketitle

\begin{center} \textbf{Abstract} \end{center}

\begin{small}
In this paper we focus on the large $n$ probability distribution
function of the largest eigenvalue in the Gaussian Orthogonal
Ensemble of $n\times n$ matrices (GOE$_n$). We prove an Edgeworth
type Theorem for the largest eigenvalue probability distribution
function of GOE$_n$. The correction terms to the limiting
probability distribution are expressed in terms of the same
Painlev$\acute{e}$ II functions appearing in the Tracy-Widom
distribution. We conclude with a brief discussion of the GSE$_{n}$
case.

\end{small}

\section{Introduction}\label{introduction}
Limiting probability distributions laws from Random Matrix Theory
have found many applications outside their initial domain of
discovery; the length of the longest increasing subsequence ( P.
Deift et al. \cite{Deif1}) properly scaled converges in distribution
to the Unitary Tracy-Widom law, the properly scaled largest
principal component of a white Wishart converges in distribution to
the Orthogonal Tracy-Widom law (I. M. Johnstone \cite{John1}). For
recent reviews we refer the reader to \cite{Deif2,Deif3, Dien1,
Johansson, Trac8}. In these applications of Tracy-Widom
distributions, one would like to control quantitatively the range of
validity of the various limit laws. One therefore needs finite $n$
correction to these limiting distributions. In a previous work
\cite{Choup1} we initiated the study of this problem for the
Gaussian Unitary Ensemble of $n$ by $n$ matrices (GUE$_{n}$).
Following this work, we will derive the analogous result for the
Gaussian Orthogonal Ensemble (GOE$_{n}$) in this paper. The
derivation of the probability distribution function of the largest
eigenvalue for Gaussian Symplectic Ensemble (GSE$_{n}$) is similar
to the GOE$_{n}$ case up to the parity of the size $n$ of matrices
in consideration. We will therefore mention at each step of the
derivation the corresponding result without much explanations except
when it is necessary. We seek a large $n$ expansion of the
probability distribution function of the largest eigenvalue in
GOE$_{n}$ and GSE$_{n}$ similar to the following Edgeworth Expansion
arising in probability in applications of the Central Limit Theorem.

If $S_n$ is a sum of  i.i.d.\  random variables $X_j$, each with
mean $\mu$ and variance $\sigma^2$,
 then the distribution $F_n$ of the normalized
random variable $(S_n - n\mu)/(\sigma \sqrt{n})$ satisfies the
Edgeworth expansion\footnote{We assume, of course,  the moments
$\bE(X_j^k)$, $k=3,\ldots, r$,  exist; and  as well, the condition
$\lim_{\vert\zeta\vert\ra\iy}\sup\vert\varphi(\zeta)\vert <\iy$
where $\varphi$ is the characteristic function of $X_j$, see
\cite{Fell2}.}
\begin{equation}\label{edgeworth}
 F_{n}(x)-\Phi(x) = \phi(x)\sum_{j=3}^{r} n^{-\frac{1}{2} j
+1}R_{j}(x)+o(n^{-\frac{1}{2}r+1})
\end{equation}
uniformly in $x$. Here $\Phi$ is the standard normal distribution
with density $\phi$, and  $R_{j}$ are  polynomials depending only on
$\bE(X_j^k)$ but not on $n$ and $r$ (or the underlying distribution
of the $X_j$).

If we view the random matrix ensembles of $n$ by $n$ matrices in
terms of the associated eigenvalues, then the Gaussian
$\beta$–-ensembles are probability spaces on $n$-tuples of random
variables $\{\lambda_{1}, \ldots , \lambda_{n}\}$ (think of them as
eigenvalues of a randomly chosen matrix from the ensemble.) The
probability density that the eigenvalues lie in an infinitesimal
intervals about the points $x_{1}, \ldots , x_{n}$ is
\begin{equation}\label{jpdfeig}
\bP_{n\beta}(x_{1},\cdots,x_{n})\; = \textrm{C}_{n\beta}\,
\textrm{exp}\left(-\frac{\beta}{2}\, \sum_{1}^{n}x_{j}^{2}\right)\,
\prod_{j<k}|x_{j}-x_{k}|^{\beta},
\end{equation}
with
\begin{equation}
-\infty < x_{i} < \infty, \quad \textrm{for} \; i=1, \cdots, n.
\end{equation}
Here $\textrm{C}_{n\beta}$ is the normalizing constant such that the
total integral over the $x_{i}$'s is one. The cases
$\beta=1,\>2,\>4$ correspond to the GOE$_{n}$, GUE$_{n}$ and
GSE$_{n}$ respectively. We denote the largest eigenvalue by
$\lambda_{Max}^{\beta}$, and
\begin{equation}\label{p.d.f}
F _{n,\beta}(t)= \bP(\lambda_{\textrm{max}}^{\beta} \leq t)
\end{equation}
the probability distribution function.\\
When $\beta=2$, the harmonic oscillator wave functions
\begin{equation*}
\varphi_{k}(x)= {1\over (2^{k}k! \sqrt{\pi})^{1/2}} \, H_k(x)\,
e^{-x^2/2}\quad k=0,\,1,\,2,\ldots
\end{equation*}
 obtained by orthonormalizing the sequence $x^{k}\,e^{-x^{2}}$ (with $H_k(x)$ the
Hermite polynomials of degree $k$) play an important role.  We also
have the Hermite kernel
\begin{equation}\label{hermite kernel}
K_{n,2}(x,y)=\sum_{k=0}^{n-1}\varphi_{k}(x)\varphi_{k}(y) =
\sqrt{\frac{n}{2}} \quad \frac{\varphi_{n}(x)\varphi_{n-1}(y)
 -  \varphi_{n}(y)\varphi_{n-1}(x)}{x-y},
\end{equation}
which is the kernel of an integral operator $K_{2}$ acting on
$L^{2}(t,\iy)$ , with resolvent
\begin{equation}\label{resolvent kernel}
R_{n}(x,y;t)\>=\> (I-K_{n})^{-1}K_{n}(x,y).
\end{equation}
The product on the right is operator multiplication.
 We have the following representation of
\eqref{jpdfeig}, (see for example, \cite{Meht1})
\begin{equation*}
\bP_{n2}(x_{1},\cdots,x_{n})\>=\> \det(K_{n}(x_{i},x_{j}))_{1\leq
i\,j\,\leq n}.
\end{equation*}
Following Tracy and Widom \cite{Trac2}, we define
\begin{equation}\label{varphi}
\varphi(x)\>=\>
\biggl(\frac{n}{2}\biggr)^{\frac{1}{4}}\varphi_{n}(x),\quad \quad
\psi(x)\>=\>
\biggl(\frac{n}{2}\biggr)^{\frac{1}{4}}\varphi_{n-1}(x),
\end{equation}
by $\varepsilon$ the integral operator with kernel
\begin{equation}\label{varepsilon}
\varepsilon_{t}(x)\>\>=\>\> \frac{1}{2}\mathrm{sgn}(x-t),
\end{equation}
$D$ the differentiation with respect to the independent variable,
\begin{equation}\label{Qn}
Q_{n,i}(x;t)\, =\, (\,(I-K_{n})^{-1}\, ,\, x^{i}\varphi_{n})
\end{equation}
and
\begin{equation}\label{Pn}
P_{n,i}(x;t)\, =\, (\,(I-K_{n})^{-1}\, ,\, x^{i}\varphi_{n-1}).
\end{equation}
We introduce the following quantities
\begin{equation}\label{qn}
q_{n,i}(t)\>=\> Q_{n,i}(t;t), \>\>\> p_{n,i}(t)\>=\> P_{n,i}(t;t)
\end{equation}
\begin{equation}\label{un}
u_{n,i}(t)\,=\, (Q_{n,i},\varphi_{n}),\quad v_{n,i}(t)\,=\,
(P_{n,i},\varphi_{n}),\quad
\end{equation}
\begin{equation}
\tilde{v}_{n,i}(t)=(Q_{n,i},\varphi_{n-1}),\quad \mathrm{and} \quad
w_{n,i}(t)\,= \, (P_{n,i},\varphi_{n-1}).
\end{equation}
 Here $(\,\cdot \,,\cdot \,)$ denotes
the inner product on $L^{2}(t,\infty)$. In our notation, the
subscript without the $n$ represents the scaled limit of that
quantity when $n$ goes to infinity, and  we dropped the second
subscript $i$ when it is zero.\\  If $\airy$ is Airy function, the
kernel $K_{n,2}(x,y)$ then scales\footnote{For the precise
definition of this scaling, see the next section} to the  Airy
kernel
\begin{equation}\label{Airy kernel}
K_{\airy}(X,Y)\>=\>
\frac{\airy(X)\,\airy'(Y)\>-\>\airy(Y)\,\airy'(X)}{X-Y}.
\end{equation}
Our conventions are as follows:
\begin{equation}\label{Q}
Q_{i}(x;s)\, =\, (\,(I-K_{\airy})^{-1}\, ,\, x^{i}\airy),\quad
Q_{0}(x;s)\,=\,Q(x;s),
\end{equation}
\begin{equation}\label{P}
P_{i}(x;s)\, =\, (\,(I-K_{\airy})^{-1}\, ,\, x^{i}\airy'),\quad
P_{0}(x;s)\,=\,P(x;s),
\end{equation}
\begin{equation}\label{q}
q_{i}(s)\>=\> Q_{i}(s;s),\quad q_{0}(s)=q(s) ,\>\>\> p_{i}(s)\>=\>
P_{i}(s;s),\quad p_{0}(s)=p(s),
\end{equation}
\begin{equation}\label{u}
u_{i}(s)\,=\, (Q_{i},\airy),\quad u_{0}(s)=u(s),\quad v_{i}(s)\,=\,
(P_{i},\airy),\quad v_{0}(s)=v(s),
\end{equation}
\begin{equation}
\tilde{v}_{i}(s)=(Q_{i},\airy'),\quad
\tilde{v}_{0}(s)=\tilde{v}(s),\quad  w_{i}(t)\,= \,
(P_{i},\airy'),\quad \textrm{and}\quad w_{0}(s)=w(s).
\end{equation}
 Here $(\,\cdot \,,\cdot \,)$  denotes
the inner product on $L^{2}(s,\infty)$  and $i=0,1,2.$ \\
We also note that $q(s)$ is the solution to the Pailev$\acute{e}$ II
equation $q''(s)=sq(s)+2q^{3}(s)$ with the boundary condition
$q(s)\sim \mathrm{Ai}(s)$ as $s \rightarrow \infty$. \\
We use the subscript $n$ for unscaled quantities only.
\begin{equation}
\mathcal{R}_{n,1}:=\int_{-\infty}^{t}R_{n}(x,t;t)\mathrm{d}x,\quad
\mathcal{P}_{n,1}:= \int_{-\infty}^{t}P_{n}(x;t)\mathrm{d}x,\quad
\mathcal{Q}_{n,1}:= \int_{-\infty}^{t}Q_{n}(x;t)\mathrm{d}x,
\end{equation}
and
\begin{equation*}
\mathcal{R}_{n,4}(t):=\int_{-\infty}^{\infty}\varepsilon_{t}(x)R_{n}(x,t;t)\mathrm{d}x,\quad
\mathcal{P}_{n,4}(t):=
\int_{-\infty}^{\infty}\varepsilon_{t}(x)P_{n}(x;t)\mathrm{d}x,
\end{equation*}
\begin{equation}
\mathcal{Q}_{n,4}(t):=
\int_{-\infty}^{\infty}\varepsilon_{t}(x)Q_{n}(x;t)\mathrm{d}x.
\end{equation}
The epsilon quantities are
\begin{equation}\label{Q epsilon}
Q_{n,\varepsilon}(x;t)\> =\> \bigl(\,(I-K_{n})^{-1}(x,y)\, ,\,
\varepsilon \varphi(y)\bigr), \quad
q_{n,\varepsilon}(t)\>=\>Q_{n,\varepsilon}(t;t)
\end{equation}

\begin{equation}\label{u epsilon}
u_{n,\varepsilon}(t)\,=\,
\bigl(Q_{n,\varepsilon}(x;t),\varphi(x)\bigr),\quad
\tilde{v}_{n,\varepsilon}(t)\,=\, \bigl(Q_{n,\varepsilon}(x;t),
\psi(x)\bigr),
\end{equation}
where $(\,\cdot \,,\cdot \,)$  denotes the inner product on
$L^{2}(t,\infty)$. \\

The GOE$_{n}$ and GSE$_{n}$ analogue of \eqref{GUE Edgeworth} in
Theorem \ref{Fn2} bellow will follow from representations (equations
(40) and (41) of \cite{Trac2}.)
\begin{equation}\label{f_{n,1}}
F_{n,1}(t)^{2}
       = F_{n,2}(t)\cdot
       \left( \bigl(1-\tilde{v}
       _{n,\varepsilon}(t)\bigr)\bigl(1-\frac{1}{2}\mathcal{R}_{n,1}(t)\bigr)-\frac{1}{2}\bigl(q_{n,\varepsilon}(t)
       -c_{\varphi}\bigr)\mathcal{P}_{n,1}(t)\right)
\end{equation}
and
\begin{equation}\label{f_{n,4}}
F_{n,4}(t/\sqrt{2})^{2}= F_{n,2}(t)\cdot\left(
       \bigl(1-\tilde{v}_{n,\varepsilon}(t)\bigr)\bigl(1+\frac{1}{2}\mathcal{R}_{n,4}(t)\bigr)+
       \frac{1}{2}q_{n,\varepsilon}(t)\,\mathcal{P}_{n,4}(t)\right).
\end{equation}
Here we first derive a large $n$-expansion of
$\mathcal{R}_{n,1},\>\mathcal{P}_{n,1},\>\mathcal{R}_{n,4},\>\mathcal{P}_{n,4},\>\tilde{v}
       _{n,\varepsilon},\>\mathrm{and}\>q_{n,\varepsilon}$ in terms
       of $p_{n}$ and $q_{n}$,  by solving the associated systems of differential
       equations. We then substitute the resulting
       expressions in \eqref{f_{n,1}} and \eqref{f_{n,4}}. We will
need the following result which gives the large $n$ expansion of
$F_{n,2}$.
\begin{thm}\cite{Choup1}\label{Fn2}
If we set
\begin{equation}\label{scale for GUE}
 t  =  (2(n+c)
)^{\frac{1}{2}}+2^{-\frac{1}{2}}n^{-\frac{1}{6}}\,s\>\>\>
\textrm{and}
\end{equation}
\begin{equation}\label{second term}
E_{c,2}:=E_{c,2}(s)=2w_{1}-3u_{2}+ (-20c^2 +3)v_{0}  +
u_{1}v_{0}-u_{0}v_{1} +u_{0}v_{0}^{2}-u_{0}^{2}w_{0}.
\end{equation}
Then as $n\ra\iy$
\begin{equation}\label{GUE Edgeworth}
F_{n,2}(t)= F_{2}(s)\left\{ 1 + c \, u_{0}(s) \, n^{-\frac{1}{3}}
-\frac{1}{20}E_{c,2}(s)\, n^{-\frac{2}{3}}\right\} + O(n^{-1})
\end{equation}
uniformly in $s$, and
\begin{equation}\label{TW}
F_{2}(s)\>=\> \lim_{n\ra\iy}F_{n,2}(t)\>=\>
\exp\left(-\int_{s}^{\iy}(x-s)q(x)^2 \,dx\right)
\end{equation}
is the Tracy-Widom distribution.
\end{thm}

\subsection{Statement of our results}

To state our main result we need the following definitions
\begin{equation}\label{alpha}
\alpha:=\alpha(s)\>=\> \int_{s}^{\iy}q(x)\,u(x)\,dx,
\end{equation}
\begin{equation}\label{mu}
\mu:=\mu(s)\>=\>\int_{s}^{\iy}q(x)dx,
\end{equation}
\begin{equation}\label{nu}
\nu:=\nu(s)\>=\>\int_{s}^{\iy}p(x)dx \>=\> \alpha(s)\>-\> q(s),
\end{equation}

\begin{equation*}
\eta:=\eta(s)=\>\frac{1}{20\sqrt{2}}\int_{s}^{\infty}\hspace{-0.1in}\bigl(6qv+3pu+2p_{2}+
2p_{1}v+2pv_{1}-2q_{2}u-2q_{1}u_{1}-2qu_{2}\bigr)(x)\,dx \,-
\end{equation*}
\begin{equation}
\frac{20c^{2}q'(s)+3p(s)}{20\sqrt{2}}
\end{equation}

\begin{equation*}
E_{c,1}(s)\,=-\frac{1}{20} E_{c,2}(s)\,e^{-\mu }\,
-\,\frac{c\,\alpha}{2 \mu ^2}+\frac{c\, p}{2 \mu }+\frac{(2c-1)\,
\nu ^2}{4 \mu ^2} +c\,u \left(c\,q\,e^{-\mu }\, - \, \frac{\nu }{2
\mu }(1-e^{-\mu } )\right) \>\>+
\end{equation*}
\begin{equation*}
e^{-2 \mu } \left(\frac{\nu\,(\nu+8c\,q) }{32 \mu } -\frac{\eta }{4
\sqrt{2}}\right)\,+ \,e^{-\mu } \left(\frac{2\sqrt{2}\,c^2
q^2-3\,\eta}{4\sqrt{2}}
\,+\,\frac{\nu^{2}-8(2c\,p+c^{2}q^{2})-4c^{2}\alpha^{2}}{32\mu}\right.
\end{equation*}
\begin{equation}\label{Ec1}
-\frac{c^{2}q^{2}}{8\mu^{2}}\,+\,\frac{2-\mu }{2 \mu ^2}
\left(c\,q\,\alpha +\frac{1}{4}\nu^{2}+(c^{2}-c)q^{2}\right)\Biggr)
-\left(4c^{2}\alpha^{2}+3c^{2}q^{2}-\nu^{2}\right)\frac{\cosh(\mu)}{8\mu^{2}}.
\end{equation}

Unfortunately, we did not find a simple representation of $\eta$ and
$E_{c,1}$. Nevertheless the quantities $\alpha(s),\> \mu(s),\>
\nu(s),\> \eta(s)$ and $E_{c,1}(s)$ are easy to compute numerically.
For $E_{c,1}(s)$ and $\eta(s)$ we only need the recurrences
relations defining $p_{i}(s)$ and $q_{i}(s)$ in term of $q(s)$ and
$q'(s)$. We find the following representation of the large $n$
probability distribution function for $\lambda_{max}$ in GOE$_{n}$.
The derivation of the analogous result for the GSE$_{n}$ follows
exactly the same steps as the one given in this paper. Here is our
main result.

\begin{thm}\label{GOE}
We set
\begin{equation}\label{scale for GOE}
t  =  (2(n+c)
)^{\frac{1}{2}}+2^{-\frac{1}{2}}n^{-\frac{1}{6}}\,s.\>\>\>
\end{equation}
Then as $n\ra\iy$
\begin{equation*}
F_{n,1}(t)^{2}\>=\>F_{2}(s)\cdot\left\{e^{-\mu(s)}\>+\>
\left[c\bigl(q(s)+u(s)\bigr) e^{-\mu(s) } -\frac{\nu(s) }{2 \mu(s)
}\bigl(1-e^{-\mu(s) }\bigr) \right]\,n^{-\frac{1}{3}} \>\> + \right.
\end{equation*}
\begin{equation}\label{Fn1}
 E_{c,1}(s) n^{-\frac{2}{3}}\>\biggr\}\>+\>O(n^{-1})
\end{equation}
uniformly for bounded $s$.
\end{thm}
Note that unlike the GUE$_{n}$ case where for $c=0$ the
$n^{-\frac{1}{3}}$ correction term vanishes as shown in equation
\eqref{GUE Edgeworth}, the  $n^{-\frac{1}{3}}$  correction term does
not vanish in the GOE$_{n}$
no matter what the fine tuning constant $c$ is.\\

 In \S 2 we reproduce the derivation of \eqref{f_{n,1}}
following Tracy and Widom in \cite{Trac2}. In \S 3 we solve the
system of equations satisfied by the various functions on the right
of \eqref{f_{n,1}} for our derivation of the Edgeworth expansion of
the probability distribution of the largest eigenvalue in the
GOE$_{n}$.

\section{Derivation of $F_{n,\beta}$}
We treat here the case $n$ even. If we set
\begin{equation}
K_{n,1}\,=\> \left(\begin{array}{cc}
                K_{n,2}+\psi\otimes \varepsilon \varphi & K_{n,2}\,D\,-\,\psi\otimes \varphi \\
                \varepsilon\,K_{n,2}-\varepsilon +\varepsilon\,\psi\otimes \varepsilon\varphi & K_{n,2}+\varepsilon\psi\otimes  \varphi \\
              \end{array}
              \right)
\end{equation}
then $F_{n,1}^{2}(t)$ is the Fredholm determinant of $K_{n,1}$ on
the set $J=(t,\iy)$. If we denote by $\chi$ the multiplication by
the function $\chi_{J}(x)$, then $F_{n,1}^{2}(t)$ is the Fredholm
determinant of the integral operator with kernel\footnote{To
simplify notation we kept the same notation for the integral
operator as well as the kernel }
\begin{equation*}
K_{n,1}\,=\> \chi_{J}(x)\,\left(\begin{array}{cc}
                K_{n,2}+\psi\otimes \varepsilon \varphi & K_{n,2}\,D\,-\,\psi\otimes \varphi \\
                \varepsilon\,K_{n,2}-\varepsilon +\varepsilon\,\psi\otimes \varepsilon\varphi & K_{n,2}+\varepsilon\psi\otimes  \varphi \\
              \end{array}
              \right)\chi_{J}(y)
\end{equation*}
\begin{equation}\label{kernel1}
=\> \chi\,\left(\begin{array}{cc}
                K_{n,2}+\psi\otimes \varepsilon \varphi & K_{n,2}\,D\,-\,\psi\otimes \varphi \\
                \varepsilon\,K_{n,2}-\varepsilon +\varepsilon\,\psi\otimes \varepsilon\varphi & K_{n,2}+\varepsilon\psi\otimes  \varphi \\
              \end{array}
              \right)\chi
\end{equation}
on $\mathbb{R}$ , see for example \cite{Trac2} equation (31). Using
the following commutators,
\begin{equation}\label{commutator}
[K_{n,2},D]\>=\> \varphi\otimes\psi +\psi \otimes \varphi,\quad
\quad [\varepsilon, K_{n,2}]\>=\> -\varepsilon \varphi \otimes
\varepsilon \psi -\varepsilon\psi \otimes \varepsilon \varphi
\end{equation}
($\psi$ and $\varphi$ appear as a consequence of the Christophel
Darboux formula applied to $K_{n,2}$,) we have
\begin{equation*}
K_{n,2}+\psi\otimes\varphi=D\,\varepsilon
K_{n,2}+D\,\varepsilon\psi\otimes\varepsilon\varphi=D(\varepsilon
K_{n,2}+\varepsilon\psi\otimes \varepsilon\varphi)=D(
K_{n,2}\,\varepsilon-\varepsilon\varphi\otimes \varepsilon\psi)
\end{equation*}
\begin{equation*}
K_{n,2}\,D-\psi\otimes\varepsilon=D
\,K_{n,2}+\varphi\otimes\psi=DK_{2}+D\varepsilon\psi\otimes\varepsilon=
D\,(K_{n,2}+\varepsilon\varphi\otimes\psi)
\end{equation*}
and
\begin{equation*}
\varepsilon\,K_{n,2}-\varepsilon +\varepsilon\,\psi\otimes
\varepsilon\varphi =K_{n,2}\,\varepsilon-\varepsilon
-\varepsilon\,\varphi\otimes \varepsilon\psi
\end{equation*}
as $D\varepsilon\>=\> I$. Our kernel is now
\begin{equation}
K_{n,1}\,=\> \chi\left(\begin{array}{cc}
                D(
K_{n,2}\,\varepsilon-\varepsilon\varphi\otimes \varepsilon\psi) & D\,(K_{n,2}+\varepsilon\varphi\otimes\psi) \\
                K_{n,2}\,\varepsilon-\varepsilon
-\varepsilon\,\varphi\otimes \varepsilon\psi & K_{n,2}+\varepsilon\psi\otimes  \varphi \\
              \end{array}
              \right)\chi
\end{equation}
\begin{equation}
=\> \left( \begin{array}{cc} \chi\,D & 0\\ 0& \chi\\
\end{array}\right)\cdot
\left(\begin{array}{cc}
                (
K_{n,2}\,\varepsilon-\varepsilon\varphi\otimes \varepsilon\psi)\chi & (K_{n,2}+\varepsilon\varphi\otimes\psi)\chi \\
                (K_{n,2}\,\varepsilon-\varepsilon
-\varepsilon\,\varphi\otimes \varepsilon\psi )\,\chi& (K_{n,2}+\varepsilon\psi\otimes  \varphi)\, \chi \\
              \end{array}
              \right)
\end{equation}
Since $K_{n,1}$ is of the form $\mathrm{AB}$, we can use the fact
that $\det(I-\mathrm{AB}) = \det(I-\mathrm{BA})$ and deduce that the
Fredholm determinant of $K_{n,1}$ is unchanged if instead we take
$K_{n,1}$ to be
\begin{equation}
\left(\begin{array}{cc}
                (
K_{n,2}\,\varepsilon-\varepsilon\varphi\otimes \varepsilon\psi)\chi & (K_{n,2}+\varepsilon\varphi\otimes\psi)\chi \\
                (K_{n,2}\,\varepsilon-\varepsilon
-\varepsilon\,\varphi\otimes \varepsilon\psi )\,\chi& (K_{n,2}+\varepsilon\psi\otimes  \varphi)\, \chi \\
              \end{array}
              \right)\cdot\left( \begin{array}{cc} \chi\,D & 0\\ 0& \chi\\
\end{array}\right)
\end{equation}
\begin{equation}
\>\>=\>\>\left(\begin{array}{cc}
                (
K_{n,2}\,\varepsilon-\varepsilon\varphi\otimes \varepsilon\psi)\chi\,D & (K_{n,2}+\varepsilon\varphi\otimes\psi)\chi \\
                (K_{n,2}\,\varepsilon-\varepsilon
-\varepsilon\,\varphi\otimes \varepsilon\psi )\,\chi\,D & (K_{n,2}+\varepsilon\psi\otimes  \varphi)\, \chi \\
              \end{array}
              \right)
\end{equation}
\begin{equation}\label{fredholm determinant}
\det(I-K_{n,1})\>=\>\det\left(\begin{array}{cc}
               I- (
K_{n,2}\,\varepsilon-\varepsilon\varphi\otimes \varepsilon\psi)\chi\,D & -(K_{n,2}+\varepsilon\varphi\otimes\psi)\chi \\
                -(K_{n,2}\,\varepsilon-\varepsilon
-\varepsilon\,\varphi\otimes \varepsilon\psi )\,\chi\,D & I-(K_{n,2}+\varepsilon\psi\otimes  \varphi)\, \chi \\
              \end{array}
              \right).
\end{equation}
Performing row and column operations on the matrix\footnote{This
does not change the determinant, for more details see \cite{Trac2}}
does not change the Fredholm determinant. We first subtract row 1
from row 2, next we add column 2 to column 1 to have the following
matrix
\begin{equation}
\left(\begin{array}{cc}
               I- (
K_{n,2}\,\varepsilon-\varepsilon\varphi\otimes
\varepsilon\psi)\chi\,D -(K_{n,2}+\varepsilon\varphi\otimes\psi)\chi
&\quad -(K_{n,2}+\varepsilon\varphi\otimes\psi)\chi \\
                \varepsilon \,\chi\,D & I \\
              \end{array}
              \right).
\end{equation}
Right-multiply column 2 by $-\varepsilon \,\chi\,D$  and add it to
column 1, and multiply row 2 by
$(K_{n,2}+\varepsilon\varphi\otimes\psi)\chi$  and add it to row 1
to have
\begin{equation}
\left(\begin{array}{cc}
               I- (
K_{n,2}\,\varepsilon-\varepsilon\varphi\otimes
\varepsilon\psi)\chi\,D -(K_{n,2}+\varepsilon\varphi\otimes\psi)\chi
+(K_{n,2}+\varepsilon\varphi\otimes\psi)\chi\,\varepsilon \,\chi\,D
&\quad 0\\   0 & \quad I \\
              \end{array}
              \right).
\end{equation}
We therefore have,
\begin{equation}\label{fredholm representation}
\det(I-K_{n,1})\>=\> \det\left(I- (
K_{n,2}\,\varepsilon-\varepsilon\varphi\otimes
\varepsilon\psi)\chi\,D
+(K_{n,2}+\varepsilon\varphi\otimes\psi)\chi\,\bigl(\varepsilon
\,\chi\,D-I\bigr)\right)
\end{equation}
\begin{equation}
=\>\det\left(I-K_{n,2}\chi
-K_{n,2}(I-\chi)\varepsilon\chi\,D-(\varepsilon\varphi\otimes
\psi)(\chi-\chi\varepsilon\chi
D)-\varepsilon\varphi\otimes\psi\varepsilon \chi D \right)
\end{equation}
We used the fact that $\varepsilon$ is antisymmetric to have
\begin{equation*}
\varepsilon\varphi\otimes
\varepsilon\psi\,\chi\,D=\varepsilon\varphi\otimes
\psi\varepsilon^{t}\,\chi\,D =-\varepsilon\varphi\otimes
\psi\varepsilon\,\chi\,D,
\end{equation*}
and if we note that $\chi$ is multiplication, then the determinant
is
\begin{equation*}
=\>\det\left(I-K_{n,2}\chi
-K_{n,2}(I-\chi)\varepsilon\chi\,D-\varepsilon\varphi\otimes
\psi(1-\chi)\varepsilon\chi D-\varepsilon\varphi\otimes\chi\psi
\right).
\end{equation*}
Next we factor out $I-K_{n,2}$ and note that
$(I-K_{n,2})^{-1}=I+R_{n,2}$, where $R_{n,2}$ was defined as the
resolvent of $K_{n,2}$, and
$(I-K_{n,2})^{-1}\varepsilon\varphi=Q_{n,\varepsilon}$. We are
interested on the determinant of the following operator
\begin{equation}\label{product determinants}
\bigl(I-K_{n,2}\chi\bigr)\bigl(I- (K_{n,2}+
R_{n,2}K_{n,2})(I-\chi)\varepsilon\chi\,D-Q_{n,\varepsilon}\otimes
\psi(1-\chi)\varepsilon\chi D-Q_{n,\varepsilon}\otimes\chi\psi
\bigr).
\end{equation}

We have a large n-expansion of $\det(I-K_{n,2}\chi)=F_{n,2}$ from
the author work in GUE$_{n}$ see \cite{Choup1} or \cite{Choup2}. We
will therefore focus our attention on the second factor of
\eqref{product determinants}. We will represent this factor in the
form $(I-\sum_{j=1}^{k}\alpha_{j}\otimes\beta_{j})$ and use the well
known formula
$\det(I-\sum_{j=1}^{k}\alpha_{j}\otimes\beta_{j})=\det\bigl(\delta_{i,j}-(\alpha_{i},\beta_{j})\bigr)_{i,j=1,\ldots
,k}$ to expand the Fredholm determinant. First we need to find a
representation of $\varepsilon\chi\,D$ as a finite rank operator. To
this end we introduce in this section the following notation,
\begin{equation*}
\varepsilon_{k}(x)=\varepsilon(x-a_{k}),\quad
R_{k}(x)=R_{n,2}(x,a_{k}),\quad \delta_{k}(x)=\delta(x-a_{k}),\quad
a_{1}=t,\quad \textrm{and}\quad a_{2}=\iy.
\end{equation*}
With the new notation $J=(t,\iy)=(a_{1},a_{2})$, and the commutator
\begin{equation*}
[\chi\>\>D]=-\delta_{1}\otimes\delta_{1}+\delta_{2}\otimes
\delta_{2},
\end{equation*}
gives
\begin{equation*}
\varepsilon
[\chi\>\>D]=-\varepsilon_{1}\otimes\delta_{1}+\varepsilon_{2}\otimes
\delta_{2}.
\end{equation*}
Next we use the identity $\varepsilon\, D=I$ to have
\begin{equation}\label{eq1}
(I-\chi)\varepsilon\,\chi\,D = (I-\chi)\varepsilon\,[\chi
\>\>D]=-(I-\chi)\varepsilon_{1}\otimes\delta_{1}+(I-\chi)\varepsilon_{2}\otimes
\delta_{2},
\end{equation}
and the representation
\begin{equation}\label{eq2}
(K_{n,2}+R_{n,2}\,K_{n,2})(I-\chi)\varepsilon\,\chi\,D=\sum_{k=1,2}(-1)^{k}(K_{n,2}+R_{n,2}\,K_{n,2})(I-\chi)
\varepsilon_{k}\otimes\delta_{k}.
\end{equation}
We substitute \eqref{eq1} and \eqref{eq2} in the second factor of
\eqref{product determinants} and have,
\begin{equation}\label{eq3}
I- \sum_{k=1,2}(-1)^{k}(K_{n,2}+
R_{n,2}K_{n,2})(I-\chi)\varepsilon_{k}\otimes\delta_{k}-\sum_{k=1,2}(-1)^{k}Q_{n,\varepsilon}\otimes
\psi\,\cdot\,(1-\chi)\varepsilon_{k}\otimes\delta_{k}
-Q_{n,\varepsilon}\otimes\chi\psi .
\end{equation}
The dot in this formula represent operator multiplication. In this
case we just multiply the kernels using the formula
$(\alpha\otimes\beta)(\gamma\otimes \delta)= (\beta,\gamma)
\alpha\otimes \delta$, to have the following form of \eqref{eq3}
\begin{equation}\label{eq4}
I- \sum_{k=1,2}(-1)^{k}(K_{n,2}+
R_{n,2}K_{n,2})(I-\chi)\varepsilon_{k}\otimes\delta_{k}-\sum_{k=1,2}(-1)^{k}\bigl(\psi,(I-\chi)\varepsilon_{k}\bigr)
Q_{n,\varepsilon}\otimes \delta_{k}
-Q_{n,\varepsilon}\otimes\chi\psi.
\end{equation}
We have
\begin{equation*}
\varepsilon_{2}=-\frac{1}{2},\quad
(I-\chi)\varepsilon_{1}=(I-\chi)\varepsilon_{2}=-\frac{1}{2},\quad
\textrm{and}\quad R_{2}=R_{n,2}(x,\iy)=0.
\end{equation*}
If we substitute these value in \eqref{eq4}, it then becomes,
\begin{equation}
I-Q_{n,\varepsilon}\otimes \chi \psi
-\frac{1}{2}\bigl[(k_{n,2}+R_{n,2}\,K_{n,2})(I-\chi)+\bigl(\psi,(I-\chi)\bigr)Q_{n,\varepsilon}\bigr]\otimes(\delta_{1}-\delta_{2}).
\end{equation}
This operator is of the desired form
\begin{equation*}
I-\sum_{k=1,2}\alpha_{k}\otimes\beta_{k}
\end{equation*}
with
\begin{equation}
\alpha_{1}=Q_{n,\varepsilon},\>\> \alpha_{2}=
\frac{1}{2}\bigl[(k_{n,2}+R_{n,2}\,K_{n,2})(I-\chi)+a_{1}Q_{n,\varepsilon}\bigr],\>\>\beta_{1}=\chi
\psi,\>\>\beta_{2}=\delta_{1}-\delta_{2},
\end{equation}
and
\begin{equation*}
a_{1}=\bigl(\psi,(I-\chi)\bigr).
\end{equation*}
The corresponding inner product are;
\begin{equation}
(\alpha_{1},\beta_{1})=\tilde{v}_{n,\varepsilon},\quad
(\alpha_{1},\beta_{2})=q_{n,\varepsilon}+c_{\varphi}
\end{equation}
with
\begin{equation}
c_{\varphi}=\varepsilon
\,\varphi(\iy)=\frac{1}{2}\int_{-\iy}^{\iy}\varphi(x)\,dx \quad
c_{\psi}=\varepsilon\,\psi(\iy)=\frac{1}{2}\int_{-\iy}^{\iy}\psi(x)\,dx,
\end{equation}
and for $n$ even
\begin{equation}\label{c varphi}
c_{\varphi}=(\pi\,n)^{1/4}2^{-3/4 -n/2}\frac{(n!)^{1/2}}{(n/2)!},
\end{equation}
and
\begin{equation}
(\alpha_{2},\beta_{1})=\frac{1}{2}(\mathcal{P}_{n,1}-a_{1}+a_{1}\tilde{v}_{n,\varepsilon}),\quad
(\alpha_{2},\beta_{2})=\frac{1}{2}(\mathcal{R}_{n,1}+a_{1}q_{n,\varepsilon}-a_{1}c_{\varphi}).
\end{equation}
The determinant of \eqref{eq4} is therefore
\begin{equation}\label{determinant}
(1-\tilde{v}_{n,\varepsilon})(1-\frac{1}{2}\mathcal{R}_{n,1})-\frac{1}{2}(q_{n,\varepsilon}-c_{\varphi})\mathcal{P}_{n,1}.
\end{equation}
In a similar way we obtain the second factor on the right side of
\eqref{f_{n,4}} for the GSE$_{n}$ case
\begin{equation}
\bigl(1-\tilde{v}_{n,\varepsilon}\bigr)\bigl(1+\frac{1}{2}\mathcal{R}_{n,4}\bigr)+
       \frac{1}{2}q_{n,\varepsilon}(t)\,\mathcal{P}_{n,4}(t).
\end{equation}

We will derive differential equations involving the various terms in
equation \eqref{determinant} and solutions in terms of $q_{n}$ and
$p_{n}$. Then used known asymptotic of $q_{n}$ and $p_{n}$ for the
derivation of a large $n$ expansion of $F_{n,1}$.

\section{Differential Equations}
\subsection{System of Differential Equations}
This section will be devoted to finding expressions of
\begin{equation*}
\mathcal{R}_{n,1}(t),\quad \mathcal{P}_{n,1}(t),\quad
q_{n,\varepsilon}(t),\quad \textrm{and}\quad
\tilde{v}_{n,\varepsilon}(t),
\end{equation*}
 and show how to obtain the corresponding quantities
\[ \mathcal{R}_{n,4}(t),\quad \textrm{and} \quad
\mathcal{P}_{n,4}(t)\] for the GSE$_{n}$.

To solve the associated system of differential equations it is
convenient to introduce the following quantities
\begin{equation*}
\mathcal{Q}_{n,1}(t),\quad  u_{n,\varepsilon}(t),
\end{equation*}
and
\begin{equation}
\rho_{n,2}(x,y) \quad \textrm{the kernel of } \>
(I-K_{n,2}\chi)^{-1}.
\end{equation}
We have
\begin{equation}
\rho_{n,2}(x,y)\>=\>\delta(x-y)\>+\>R_{n,2}(x,y),
\end{equation}
and
\begin{equation}
\frac{d}{dt}\mathcal{R}_{n,1}(t)=\frac{d}{dt}\int_{-\iy}^{t}R_{n,2}(x,t)\,dx=R_{n,2}(t,t)+
\int_{-\iy}^{t}\frac{d}{dt}R_{n,2}(x,t)\,dx.
\end{equation}
Formula (45) of \cite{Trac2} gives
\begin{equation}
\frac{d}{dt}R_{n,2}(x,t)=-\frac{d}{dx}R_{n,2}(x,t)-p_{n}(t)Q_{n}(x;t)-q_{n}P_{n}(x;t),
\end{equation}
and we find that
\begin{equation}
\mathcal{R}_{n,1}^{'}(t)\>=\>\frac{d}{dt}\mathcal{R}_{n,1}(t)=-p_{n}(t)\,\mathcal{Q}_{n,1}(t)-q_{n}(t)\,\mathcal{P}_{n,1}(t).
\end{equation}
We also have
\begin{equation}
\mathcal{Q}_{n,1}^{'}(t)\>=\>\frac{d}{dt}\mathcal{Q}_{n,1}(t)\>=\>\frac{d}{dt}\int_{-\iy}^{t}Q_{n}(x;t)\,dx
\>=\> q_{n}(t)(1-\mathcal{R}_{n,1}(t)),
\end{equation}
and
\begin{equation}
\mathcal{P}_{n,1}^{'}(t)\>=\>\frac{d}{dt}\mathcal{P}_{n,1}(t)\>=\>\frac{d}{dt}\int_{-\iy}^{t}P_{n}(x;t)\,dx
\>=\>p_{n}(t)(1-\mathcal{R}_{n,1}(t)),
\end{equation}
where we used
\begin{equation}
\frac{d}{dt}Q_{n}(x;t)=-q_{n}(t)\,R_{n,2}(x,t),\quad
\textrm{and}\quad \frac{d}{dt}P_{n}(x;t)=-p_{n}(t)\,R_{n,2}(x,t).
\end{equation}
The other derivatives are
\begin{equation}
\frac{d}{dt}u_{n,\varepsilon}(t)=\frac{d}{dt}\int_{t}^{\iy}Q_{n}(x,t)\varepsilon\,\varphi(x)\,dx=
-q_{n}(t)\varepsilon\varphi(t)
+\int_{t}^{\iy}\frac{d}{dt}Q_{n}(x,t)\varepsilon\,\varphi(x)\,dx
\end{equation}
\begin{equation}
=-q_{n}(t)\bigl(\varepsilon\varphi(t)
+\int_{t}^{\iy}\frac{d}{dt}R_{n,2}(x,t)\varepsilon\,\varphi(x)\,dx\bigr)=
-q_{n}(t)\int_{t}^{\iy}\rho_{n,2}(x,t)\varepsilon\,\varphi(x)\,dx,
\end{equation}
and therefore
\begin{equation}
u_{n,\varepsilon}^{'}(t)\>\>=\>\>-q_{n}(t)\,q_{n,\varepsilon}(t).
\end{equation}
Similarly
\begin{equation}
\frac{d}{dt}\tilde{v}_{n,\varepsilon}(t)=\frac{d}{dt}\int_{t}^{\iy}P_{n}(x,t)\varepsilon\,\varphi(x)\,dx=
-p_{n}(t)\varepsilon\varphi(t)
+\int_{t}^{\iy}\frac{d}{dt}P_{n}(x,t)\varepsilon\,\varphi(x)\,dx,
\end{equation}
or
\begin{equation}
\tilde{v}_{n,\varepsilon}^{'}(t)\>\>=\>\>-p_{n}(t)\,q_{n,\varepsilon}(t).
\end{equation}
The last of these is
\begin{equation}
\frac{d}{dt}q_{n,\varepsilon}(t)=\frac{d}{dt}\int
\rho_{n,2}(t,y)\,\varepsilon\varphi(y)\,dy
\end{equation}
\begin{equation}
=-\int \frac{\partial}{\partial
y}\rho_{n,2}(t,y)\,\varepsilon\varphi(y)dy -q_{n}(t)\bigl(\chi
P_{n}(y;t),\varepsilon\varphi(y)\bigr)-p_{n}(t)\bigl(\chi
Q_{n}(y;t),\varepsilon\varphi(y)\bigr).
\end{equation}
Integration by parts together with the boundary conditions and
$D\varepsilon=I$ gives
\begin{equation}
-\int \frac{\partial}{\partial
y}\rho_{n,2}(t,y)\,\varepsilon\varphi(y)dy= \int
\rho_{n,2}(t,y)\,\chi\varepsilon\varphi(y)\,dy=q_{n}(t),
\end{equation}
which in turn gives
\begin{equation}
q_{n,\varepsilon}^{'}(t) =
q_{n}(t)-\tilde{v}_{n,\varepsilon}(t)q_{n}(t)-u_{n,\varepsilon}(t)q_{n}(t).
\end{equation}
The boundary conditions at $t=\iy$ for these function are,
\begin{equation}
\mathcal{R}_{n,1}(\iy)=0,\quad \mathcal{Q}_{n,1}(\iy)=2c_{\varphi}
\quad \textrm{and} \quad \mathcal{P}_{n,1}(\iy)=2c_{\psi}=0\quad
\textrm{ as }n \textrm{ is even.}
\end{equation}
and
\begin{equation}
\tilde{v}_{n,\varepsilon}(\iy)=0,\quad
u_{n,\varepsilon}(\iy)=0,\quad \textrm{and}\quad
q_{n,\varepsilon}(\iy)=c_{\varphi}.
\end{equation}
The associated systems of equations are;
\begin{equation}\label{system1}
\left\{
  \begin{array}{lll}
   q_{n,\varepsilon}^{'}(t)  &=&q_{n}(t)\bigl(1-\tilde{v}_{n,\varepsilon}(t)\bigr)-p_{n}(t)\,u_{n,\varepsilon}(t); \\
   \bigl(1- \tilde{v}_{n,\varepsilon}\bigr)^{'}(t)&= & p_{n}(t)\,q_{n,\varepsilon}(t); \\
    u_{n,\varepsilon}^{'}(t)&= & -q_{n}(t)\,q_{n,\varepsilon}(t).
  \end{array}
\right.
\end{equation}
and
\begin{equation}\label{system2}
\left\{
  \begin{array}{lll}
 \bigl(1-  \mathcal{R}_{n,1}\bigr)^{'}(t)&=& p_{n}(t)\,\mathcal{Q}_{n,1}(t)+q_{n}(t)\,\mathcal{P}_{n,1}(t); \\
   \mathcal{Q}_{n,1}^{'}(t)&=& q_{n}(t)(1-\mathcal{R}_{n,1}(t)); \\
    \mathcal{P}_{n,1}^{'}(t)&=& p_{n}(t)(1-\mathcal{R}_{n,1}(t)).
  \end{array}
\right.
\end{equation}

\subsection{Asymptotic solutions}
We will define in this subsection only
\begin{equation}
V_{n,\varepsilon}(t)=1-\tilde{v}_{n,\varepsilon}(t),\quad
\textrm{and} \quad
\tilde{\mathcal{R}}_{n,1}(t)=1-\mathcal{R}_{n,1}(t).
\end{equation}
With this notation, system \eqref{system1} is
\begin{equation}
\frac{d}{dt}\left(
  \begin{array}{c}
    u_{n,\varepsilon}(t) \\
    V_{n,\varepsilon}(t) \\
    q_{n,\varepsilon}(t) \\
  \end{array}
\right) =\left(
               \begin{array}{ccc}
                 0 & 0 & -q_{n}(t) \\
                 0 & 0 & p_{n}(t) \\
                 -p_{n}(t) & q_{n}(t) & 0 \\
               \end{array}
             \right)\,\cdot\,
\left(
  \begin{array}{c}
    u_{n,\varepsilon}(t) \\
    V_{n,\varepsilon}(t) \\
    q_{n,\varepsilon}(t) \\
  \end{array}
\right)
\end{equation}
and \eqref{system2} is
\begin{equation}
\frac{d}{dt}\left(
  \begin{array}{c}
    \mathcal{Q}_{n,1}(t) \\
    \mathcal{P}_{n,1}(t) \\
    \tilde{\mathcal{R}}_{n,1}(t) \\
  \end{array}
\right) =\left(
               \begin{array}{ccc}
                 0 & 0 & q_{n}(t) \\
                 0 & 0 & p_{n}(t) \\
                 p_{n}(t) & q_{n}(t) & 0 \\
               \end{array}
             \right)\,\cdot\,
\left(
  \begin{array}{c}
    \mathcal{Q}_{n,1}(t) \\
    \mathcal{P}_{n,1}(t) \\
    \tilde{\mathcal{R}}_{n,1}(t) \\
  \end{array}
\right).
\end{equation}
If we let
\begin{equation*}
X^{t}=\bigl(u_{n,\varepsilon}(t), V_{n,\varepsilon}(t),
q_{n,\varepsilon}(t) \bigr)\quad \textrm{and}\quad
Y^{t}=\bigl(\mathcal{Q}_{n,1}(t),
    \mathcal{P}_{n,1}(t),
    \tilde{\mathcal{R}}_{n,1}(t)\bigr),
\end{equation*}
then \eqref{system1} and \eqref{system2} have the following
representations
\begin{equation}\label{eq5}
X^{'}(t)\>=\>A(t)\,X(t)\quad \textrm{and} \quad
Y^{'}(t)\>=\>B(t)\,Y(t)
\end{equation}
with
\begin{equation}
X^{t}(\iy)\>=\> (0,1,c_{\varphi}),\quad \textrm{and}\quad
Y^{t}(\iy)=(2c_{\varphi},0,1).
\end{equation}
We note that $A(t)$ is continuous for $t$ bounded away from $-\iy$.
We need to show that our matrix $A(t)$ is bounded as an operator on
$L^{1}(t,\iy)\otimes L^{1}(t,\iy)\otimes L^{1}(t,\iy)$ for this end
we will use the Max norm. The entries of $A(t)$ are $\pm q_{n}(x)$
and $\pm p_{n}(x)$.
\begin{equation}\label{bound on qn}
\int_{t}^{\iy}|q_{n}(x)|\,dx=\frac{1}{\sqrt{2}}\int_{s}^{\iy}\bigl|q(x)+f(x)n^{-\frac{1}{3}}\bigr|dx=M_{1}
\quad \textrm{with}\quad M_{1} \,<\, \iy
\end{equation}
We made use of the following change of variables together with
formula\footnote{We set the constant $c$ in (2.29) to zero, and use
the known asymptotic of q(x) at infinity to deduce the existance of
the integral.}(2.29) of \cite{Choup2}
\begin{equation}
x=\sqrt{2n}+\frac{X}{\sqrt{2}n^{\frac{1}{6}}},\quad \textrm{and}
\quad t=\sqrt{2n}+\frac{s}{\sqrt{2}n^{\frac{1}{6}}},
\end{equation}
the fact that the asymptotics for $p$ at infinity can be obtained
from the following representation $p=q^{'}+uq$ where $u(\iy)=0$ or
that $u(x)$ is bounded for $x$ away from minus infinity. We also
assumed without lost of generalities for this section that $p(x)\sim
\mathrm{Const}\cdot x^{1/4}e^{-\frac{2}{3}x^{\frac{3}{2}}}\quad
\textrm{as} \quad x\rightarrow \iy$ which is a consequence of the
asymptotics $q(x)\sim
\frac{1}{2\sqrt{\pi}x^{1/4}}e^{-\frac{2}{3}x^{\frac{3}{2}}}\quad
\textrm{as} \quad x\rightarrow \iy$. We also remarked that the
scaled value of $q_{n}$ in (2.29) of \cite{Choup2} is represented in
terms of finite combinations of bounded functions (the $u_{i}$'s,
$v_{i}$'s, $w_{i}$'s ,) $x^{i},\quad i=0,1,2$, with $p$ and $q$. So
again we assumed that the scaled value of $q_{n}(x)$  was of order
$n^{\frac{1}{6}}x^{2}e^{-\frac{2}{3}x^{\frac{3}{2}}}\quad
\textrm{as} \quad x\rightarrow \iy$. A similar argument hold for
$p_{n}$, here we use formula (2.30) of \cite{Choup2} instead.
\begin{equation}\label{bound on pn}
\int_{t}^{\iy}|p_{n}(x)|\,dx=\frac{1}{\sqrt{2}}\int_{s}^{\iy}\bigl|q(x)+g(x)n^{-\frac{1}{3}}\bigr|dx=M_{2}
\quad \textrm{with}\quad M_{2} \,<\, \iy .
\end{equation}
Note that $||A||_{Max}$ and $||B||_{Max}$  are at most
$2M_{1}+2M_{2}$. The fundamental local existence of solution for
linear Ordinary Differential Equation says that equations
\eqref{eq5} have solutions in $(a,\iy)$ with $a$ bounded away from
infinity given by
\begin{equation}\label{fundamental solution}
X(t)= \exp\left(-\int_{t}^{\iy}A(x)dx\right)\cdot X(\iy),\quad
\textrm{and} \quad Y(t)= \exp\left(-\int_{t}^{\iy}B(x)dx\right)\cdot
Y(\iy).
\end{equation}
This solution is convenient for the large $n$ expansion of the
probability distribution since it allows us to give a series
expansion of the solutions of our solution in terms of $q_{n}$ and
$p_{n}$. The other advantage is the built in symmetries in matrices
$A$ and $B$. These symmetries make the computation of the matrix
exponential very easy. We will start with the first system $ X(t)=
\exp\left(-\int_{t}^{\iy}A(x)dx\right)\cdot X(\iy)$. We set
\begin{displaymath}
\exp\left(-\int_{t}^{\iy}A(x)\,dx\right)= \exp\left\{\left(
               \begin{array}{ccc}
                 0 & 0 & \int_{t}^{\iy}q_{n}(x)\,dx \\\\
                 0 & 0 & -\int_{t}^{\iy}p_{n}(x)\,dx \\ \\
                 \int_{t}^{\iy}p_{n}(x) dx& -\int_{t}^{\iy}q_{n}(x)\,dx& 0 \\
               \end{array}
             \right)\right\}
\end{displaymath}
\begin{equation}
=\exp(M).
\end{equation}

$M$ is of the form
\begin{equation}
M=\left(
               \begin{array}{ccc}
                 0 & 0 & a \\
                 0 & 0 & -b \\
                 b & -a & 0 \\
               \end{array}
             \right)
\end{equation}
with
\begin{displaymath}
\exp(M)=\left(
               \begin{array}{ccc}
                 1+\sum_{k\geq 1}\frac{2^{k-1}a^{k}b^{k}}{(2k)!} & -\sum_{k\geq 1}\frac{2^{k-1}a^{k+1}b^{k-1}}{(2k)!} &
\sum_{k\geq 0}\frac{2^{k}a^{k+1}b^{k}}{(2k+1)!}  \\ \\
                 -\sum_{k\geq 1}\frac{2^{k-1}a^{k-1}b^{k+1}}{(2k)!}  & 1+\sum_{k\geq 1}\frac{2^{k-1}a^{k}b^{k}}{(2k)!}  &
-\sum_{k\geq 0}\frac{2^{k}a^{k}b^{k+1}}{(2k+1)!}  \\ \\
                 \sum_{k\geq 0}\frac{2^{k}a^{k}b^{k+1}}{(2k+1)!}  & -\sum_{k\geq 0}\frac{2^{k}a^{k+1}b^{k}}{(2k+1)!}  &
1+\sum_{k\geq 1}\frac{2^{k}a^{k}b^{k}}{(2k)!}  \\
               \end{array}
             \right)
\end{displaymath}
\begin{equation}\label{series solution}
=\left(\begin{array}{ccc} \exp(M)_{11} &\exp(M)_{12}&\exp(M)_{13}\\
\exp(M)_{21}&\exp(M)_{22}&\exp(M)_{23}\\
\exp(M)_{31}&\exp(M)_{32}&\exp(M)_{33}\\
\end{array}
\right)
\end{equation}
We have
\begin{equation}\label{solution u_{n}}
u_{n,\varepsilon}(t)\>\>=\>\>
\exp(M)_{12}\>\>+\>\>c_{\varphi}\,\exp(M)_{13},
\end{equation}
\begin{equation}\label{solution V_{n}}
V_{n,\varepsilon}(t)\>\>=\>\>
\exp(M)_{22}\>\>+\>\>c_{\varphi}\,\exp(M)_{23},
\end{equation}
and
\begin{equation}\label{solution q_{n}}
q_{n,\varepsilon}(t)\>\>=\>\>
\exp(M)_{32}\>\>+\>\>c_{\varphi}\,\exp(M)_{33}.
\end{equation}

\subsubsection{Scaling}
At this point we scale the functions involved in \eqref{series
solution} in terms of $n$ at the point corresponding to the expected
value of the largest eigenvalue. If we set
\begin{equation}\label{scaling}
t=\tau(s)=\sqrt{2(n+c)}+\frac{s}{2^{\frac{1}{2}}n^{\frac{1}{6}}},
\end{equation}
then equations (2.29) and (2.30) of \cite{Choup2} are
\begin{equation*}
q_{n}(\tau(s))\,=\,Q_{n}(\tau(s);\tau(s))=n^{\frac{1}{6}}\left(
q(s)+ \left[\frac{2c-1}{2}p(s)-c
q(s)u(s)\right]n^{\frac{1}{3}}\right.
\end{equation*}
\begin{equation*}
+\left[(10c^{2}-10c+\frac{3}{2})q_{1}(s)+p_{2}(s) +
(-30c^{2}+10c+\frac{3}{2})q(s) v(s) \right.
\end{equation*}
\begin{equation*}
+  p_{1}(s) v(s) +p(s) v_{1}(s)-q_{2}(s) u(s)-q_{1}(s) u_{1}(s)-q(s)
u_{2}(s)
\end{equation*}
\begin{equation}\label{q_{n} for large n}
+ \left.\left.(-10c^{2}+\frac{3}{2})p(s) u(s) +20c^{2}q(s) u^{2}(s)
\right]\frac{n^{-\frac{2}{3}}}{20} +O(n^{-1})e_{q}(s)\right),
\end{equation}
and
\begin{equation*}
p_{n}(\tau(s))\>=\>P_{n}(\tau(s);\tau(s))=n^{\frac{1}{6}}\left(
q(s)+ \left[\frac{2c+1}{2}p(s)-c
q(s)u(s)\right]n^{\frac{1}{3}}\right.
\end{equation*}
\begin{equation*}
+\left[(10c^{2}+10c+\frac{3}{2})q_{1}(s)+p_{2}(s) +
(-30c^{2}-10c+\frac{3}{2})q(s) v(s) \right.
\end{equation*}
\begin{equation*}
+  p_{1}(s) v(s) +p(s) v_{1}(s)-q_{2}(s) u(s)-q_{1}(s) u_{1}(s)-q(s)
u_{2}(s)
\end{equation*}
\begin{equation}\label{p_{n} for large n}
+ \left.\left.(-10c^{2}+\frac{3}{2})p(s) u(s) +20c^{2}q(s) u^{2}(s)
\right]\frac{n^{-\frac{2}{3}}}{20} +O(n^{-1})e_{p}(s)\right).
\end{equation}
If we change the variable in $a$ and $b$ by setting\footnote{We use
the same letter in both sides in the change here to simplify
notation.} $x:=\tau(x)$, we obtain
\begin{equation*}
a=\int_{t}^{\iy}q_{n}(x)\,dx\,=\frac{1}{\sqrt{2}}\int_{s}^{\iy}\left(
q(x)+ \left[\frac{2c-1}{2}p(x)-c
q(x)u(x)\right]n^{\frac{1}{3}}\right.
\end{equation*}
\begin{equation*}
+\left[(10c^{2}-10c+\frac{3}{2})q_{1}(x)+p_{2}(x) +
(-30c^{2}+10c+\frac{3}{2})q(x) v(x) \right.
\end{equation*}
\begin{equation*}
+  p_{1}(x) v(x) +p(x) v_{1}(x)-q_{2}(x) u(x)-q_{1}(x) u_{1}(x)-q(x)
u_{2}(x)
\end{equation*}
\begin{equation}\label{a_{n}}
+ \left.\left.(-10c^{2}+\frac{3}{2})p(x) u(x) +20c^{2}q(x) u^{2}(x)
\right]\frac{n^{-\frac{2}{3}}}{20} +O(n^{-1})e_{q}(x)\right)\,dx,
\end{equation}
\begin{equation*}
=a_{0}(s)+a_{1}(s)n^{-1/3}+a_{2}(s)n^{-2/3} +a_{3}(s)n^{-1}
\end{equation*}
 and
\begin{equation*}
b=\int_{t}^{\iy}p_{n}(x)\,dx\>=\frac{1}{\sqrt{2}}\int_{s}^{\iy}\left(
q(x)+ \left[\frac{2c+1}{2}p(x)-c
q(x)u(x)\right]n^{\frac{1}{3}}\right.
\end{equation*}
\begin{equation*}
+\left[(10c^{2}+10c+\frac{3}{2})q_{1}(x)+p_{2}(x) +
(-30c^{2}-10c+\frac{3}{2})q(x) v(x) \right.
\end{equation*}
\begin{equation*}
+  p_{1}(x) v(x) +p(x) v_{1}(x)-q_{2}(x) u(x)-q_{1}(x) u_{1}(x)-q(x)
u_{2}(x)
\end{equation*}
\begin{equation}\label{b_{n}}
+ \left.\left.(-10c^{2}+\frac{3}{2})p(x) u(x) +20c^{2}q(x) u^{2}(x)
\right]\frac{n^{-\frac{2}{3}}}{20} +O(n^{-1})e_{p}(x)\right)\,dx
\end{equation}
\begin{equation*}
=a_{0}(s)+b_{1}(s)n^{-1/3}+b_{2}(s)n^{-2/3} +b_{3}(s)n^{-1}.
\end{equation*}
We next focus on the following expression
\begin{equation*}
a^{k}b^{k}= (ab)^{k}= (a_{0}^{2}+a_{0}(a_{1}+b_{1})n^{-\frac{1}{3}}
+ (a_{0}(b_{2}+a_{2})+a_{1}b_{1})n^{-\frac{2}{3}} +D n^{-1})^{k}.
\end{equation*}
An expansion of this expression is
\begin{equation*}
a^{k}b^{k}= (ab)^{k}=
a_{0}^{2k}+ka_{0}^{2k-2}\left(a_{0}(a_{1}+b_{1})n^{-\frac{1}{3}} +
(a_{0}(b_{2}+a_{2})+a_{1}b_{1})n^{-\frac{2}{3}} +D n^{-1}\right)+
\end{equation*}
\begin{equation*}
\frac{k(k-1)}{2}a_{0}^{2k-4}\left(a_{0}(a_{1}+b_{1})n^{-\frac{1}{3}}
+ (a_{0}(b_{2}+a_{2})+a_{1}b_{1})n^{-\frac{2}{3}} +D
n^{-1}\right)^{2}+
\end{equation*}
\begin{equation*}
\sum_{i=3}^{k}a_{0}^{2k-2i}\binom{k}{i}\left(a_{0}(a_{1}+b_{1})n^{-\frac{1}{3}}
+ (a_{0}(b_{2}+a_{2})+a_{1}b_{1})n^{-\frac{2}{3}} +D
n^{-1}\right)^{i}.
\end{equation*}
If we note that for $i\geq 3$
\begin{equation*}
\left(a_{0}(a_{1}+b_{1})n^{-\frac{1}{3}} +
(a_{0}(b_{2}+a_{2})+a_{1}b_{1})n^{-\frac{2}{3}} +D
n^{-1}\right)^{i}=O(n^{-1}),
\end{equation*}
then the sum in this last term can be represented as
\begin{equation*}
\sum_{i=3}^{k}a_{0}^{2k-2i}\binom{k}{i}\,O(n^{-1})
=\left(\sum_{i=3}^{k}a_{0}^{2k-2i}\binom{k}{i} +a_{0}^{2k}
+ka_{0}^{2k-2} +\frac{k(k-1)}{2}a_{0}^{2k-4}-
a_{0}^{2k}-ka_{0}^{2k-2}\right.
\end{equation*}
\begin{equation*}
\left.-\frac{k(k-1)}{2}a_{0}^{2k-4} \right)O(n^{-1})=\left(
(a_{0}^{2}+1)^{k}-
a_{0}^{2k}-ka_{0}^{2k-2}-\frac{k(k-1)}{2}a_{0}^{2k-4}
\right)O(n^{-1}).
\end{equation*}
We have
\begin{equation}
\frac{ka_{0}^{2k-2}}{(2k)!}\left(a_{0}(a_{1}+b_{1})n^{-\frac{1}{3}}
+ (a_{0}(b_{2}+a_{2})+a_{1}b_{1})n^{-\frac{2}{3}} +D n^{-1}\right)=
\end{equation}
\begin{equation*}
\frac{a_{0}^{2k-1}}{2(2k-1)!}\left((a_{1}+b_{1})n^{-\frac{1}{3}}+(a_{2}+b_{2})n^{-\frac{2}{3}}\right)
+\frac{a_{0}^{2k-2}}{2(2k-1)!}a_{1}b_{1}n^{-\frac{2}{3}}
+\frac{a_{0}^{2k-2}}{2(2k-1)!} D n^{-1},
\end{equation*}
and for $k\geq 2$
\begin{equation}
\frac{k(k-1)}{2(2k)!}a_{0}^{2k-4}\left(a_{0}(a_{1}+b_{1})n^{-\frac{1}{3}}
+ (a_{0}(b_{2}+a_{2})+a_{1}b_{1})n^{-\frac{2}{3}} +D
n^{-1}\right)^{2}=
\end{equation}
\begin{equation*}
\left(\frac{a_{0}^{2k-2}}{8(2k-2)!}-
\frac{a_{0}^{2k-2}}{8(2k-1)!}\right)\left((a_{1}+b_{1})^{2}n^{-\frac{2}{3}}+O(n^{-1})\right).
\end{equation*}
We have at this stage,
\begin{equation*}
1+\sum_{k\geq1}\frac{2^{k-1}a^{k}b^{k}}{(2k)!}=1+\sum_{k\geq
1}\frac{2^{k-1}a_{0}^{2k}}{(2k)!}\,+\,\sum_{k\geq
1}\frac{2^{k-1}a_{0}^{2k-1}}{2(2k-1)!}(a_{1}+b_{1})n^{-\frac{1}{3}}+\left[
\sum_{k\geq
1}\frac{2^{k-1}a_{0}^{2k-1}}{2(2k-1)!}(a_{2}+b_{2})\right.
\end{equation*}
\begin{equation*}
\left.+\sum_{k\geq 1}\frac{2^{k-1}a_{0}^{2k-2}}{2(2k-1)!}a_{1}b_{1}
+\left(\sum_{k\geq
2}\frac{2^{k-1}a_{0}^{2k-2}}{8(2k-2)!}-\sum_{k\geq
2}\frac{2^{k-1}a_{0}^{2k-2}}{8(2k-1)!}\right)(a_{1}+b_{1})^{2}\right]n^{-\frac{2}{3}}+
\end{equation*}
\begin{equation}
\sum_{k\geq 1}\frac{2^{k-1}(a_{0}^{2}+1)^{k}}{(2k)!}\,O(n^{-1})
\end{equation}
\begin{equation*}
=1+\frac{1}{2}\sum_{k\geq 1}\frac{(\sqrt{2}a_{0})^{2k}}{(2k)!}\,+\,
\frac{(a_{1}+b_{1})}{2\sqrt{2}}\sum_{k\geq
1}\frac{(\sqrt{2}a_{0})^{2k-1}}{(2k-1)!}\,n^{-\frac{1}{3}}\,+\,
\left[ \frac{(a_{2}+b_{2})}{2\sqrt{2}}\sum_{k\geq
1}\frac{(\sqrt{2}a_{0})^{2k-1}}{(2k-1)!}\right. +
\end{equation*}
\begin{equation*}
\left.\frac{a_{1}b_{1}}{2\sqrt{2}a_{0}}\sum_{k\geq
1}\frac{(\sqrt{2}a_{0})^{2k-1}}{(2k-1)!}+\left(\sum_{k\geq
1}\frac{(\sqrt{2}a_{0})^{2k}}{8(2k)!}-\sum_{k\geq
1}\frac{(\sqrt{2}a_{0})^{2k}}{8(2k+1)!}\right)(a_{1}+b_{1})^{2}\right]
n^{-\frac{2}{3}}+
\end{equation*}
\begin{equation}
\sum_{k\geq 1}\frac{(2a_{0}^{2}+2)^{k}}{(2k)!}\,O(n^{-1})
\end{equation}
\begin{equation*}
=\frac{1}{2}\left(1+\cosh(\sqrt{2}a_{0})\right)+\frac{(a_{1}+b_{1})}{2\sqrt{2}}
\sinh(\sqrt{2}a_{0})\,n^{-\frac{1}{3}}\,+\,\left[
\frac{(a_{2}+b_{2})}{2\sqrt{2}}\sinh(\sqrt{2}a_{0})\right. +
\end{equation*}
\begin{equation*}
\left.\frac{a_{1}b_{1}}{2\sqrt{2}a_{0}}\sinh(\sqrt{2}a_{0})+\frac{1}{8}\left(\cosh(\sqrt{2}a_{0})
-\frac{1}{\sqrt{2}a_{0}}\sinh(\sqrt{2}a_{0})\right)(a_{1}+b_{1})^{2}\right]
n^{-\frac{2}{3}}+
\end{equation*}
\begin{equation}\label{matrix 11}
\cosh(2a_{0}^{2}+1)\,O(n^{-1})\>\>=\>\>\exp(M)_{11}\>\>=\>\>\exp(M)_{22}
\end{equation}
A similar argument gives
\begin{equation*}
-\sum_{k\geq
1}\frac{2^{k-1}a^{k+1}b^{k-1}}{(2k)!}=\frac{1}{2}\left(1-\cosh(\sqrt{2}a_{0})\right)+\biggl[\frac{a_{1}-b_{1}}{2a_{0}}
\bigl(1-\cosh(\sqrt{2}a_{0})\bigr)\> -
\end{equation*}
\begin{equation*}
\left.
\frac{a_{1}+b_{1}}{2\sqrt{2}}\sinh(\sqrt{2}a_{0})\right]\,n^{-\frac{1}{3}}
+\left[\frac{b_{1}(b_{1}-a_{1})}{2a_{0}^{2}}+\frac{a_{2}-b_{2}}{2a_{0}}+\left(\frac{b_{1}(a_{1}-b_{1})}{2a_{0}^{2}}
-\frac{(a_{1}+b_{1})^{2}}{8}\>+\right.\right.
\end{equation*}
\begin{equation*}
\left.\left.\frac{b_{2}-a_{2}}{2a_{0}}\right)\cosh(\sqrt{2}a_{0})
+\left(\frac{5\sqrt{2}b_{1}^{2}}{16a_{0}}-\frac{3\sqrt{2}a_{1}^{2}}{16a_{0}}-\frac{\sqrt{2}a_{1}b_{1}}{8a_{0}}
-\frac{\sqrt{2}(a_{2}+b_{2})}{4}\right)\sinh(\sqrt{2}a_{0})\right]n^{-\frac{2}{3}}
\end{equation*}
\begin{equation}\label{matrix 12}
+\cosh(\sqrt{2}a_{0})\,O(n^{-1})\>\>=\>\>\exp(M)_{12},
\end{equation}
if we interchange $a$ and $b$, then
\begin{equation*}
-\sum_{k\geq
1}\frac{2^{k-1}a^{k-1}b^{k+1}}{(2k)!}=\frac{1}{2}\left(1-\cosh(\sqrt{2}a_{0})\right)+\biggl[\frac{b_{1}-a_{1}}{2a_{0}}
\bigl(1-\cosh(\sqrt{2}a_{0})\bigr)\> -
\end{equation*}
\begin{equation*}
\left.
\frac{a_{1}+b_{1}}{2\sqrt{2}}\sinh(\sqrt{2}a_{0})\right]\,n^{-\frac{1}{3}}
+\left[\frac{a_{1}(a_{1}-b_{1})}{2a_{0}^{2}}+\frac{b_{2}-a_{2}}{2a_{0}}+\left(\frac{a_{1}(b_{1}-a_{1})}{2a_{0}^{2}}
-\frac{(a_{1}+b_{1})^{2}}{8}\>+\right.\right.
\end{equation*}
\begin{equation*}
\left.\left.\frac{a_{2}-b_{2}}{2a_{0}}\right)\cosh(\sqrt{2}a_{0})
+\left(\frac{5\sqrt{2}a_{1}^{2}}{16a_{0}}-\frac{3\sqrt{2}b_{1}^{2}}{16a_{0}}-\frac{\sqrt{2}a_{1}b_{1}}{8a_{0}}
-\frac{\sqrt{2}(a_{2}+b_{2})}{4}\right)\sinh(\sqrt{2}a_{0})\right]n^{-\frac{2}{3}}
\end{equation*}
\begin{equation}\label{matrix 21}
+\cosh(\sqrt{2}a_{0})\,O(n^{-1})\>\>=\>\>\exp(M)_{21},
\end{equation}
We also have
\begin{equation*}
\sum_{k\geq
0}\frac{2^{k}a^{k+1}b^{k}}{(2k+1)!}=\frac{1}{\sqrt{2}}\sinh(\sqrt{2}a_{0})+\left[\frac{a_{1}+b_{1}}{2}\cosh(\sqrt{2}a_{0})
+\frac{a_{1}-b_{1}}{2\sqrt{2}a_{0}}\sinh(\sqrt{2}a_{0})\right]n^{-\frac{1}{3}}\>+
\end{equation*}
\begin{equation*}
\left[\left(\frac{(a_{1}+b_{1})^{2}}{8a_{0}}-\frac{b_{1}^{2}}{2a_{0}}+\frac{a_{2}+b_{2}}{2}\right)\cosh(\sqrt{2}a_{0})
+\left(\frac{(a_{1}+b_{1})^{2}}{4\sqrt{2}}-\frac{(a_{1}+b_{1})^{2}}{8\sqrt{2}a_{0}^{2}}+\frac{b_{1}^{2}}{2\sqrt{2}a_{0}^{2}}\right.
\right.
\end{equation*}
\begin{equation}\label{matrix 13}
\left.\left.
+\>\>\frac{a_{2}-b_{2}}{2\sqrt{2}a_{0}}\right)\sinh(\sqrt{2}a_{0})\right]n^{-\frac{2}{3}}+\sinh(\sqrt{2}a_{0})O(n^{-1})
\>\>=\>\>\exp(M)_{13}\>\>=\>\>-\exp(M)_{32},
\end{equation}
and if we interchange $a$ and $b$ in this last formula,
\begin{equation*}
\sum_{k\geq
0}\frac{2^{k}a^{k}b^{k+1}}{(2k+1)!}=\frac{1}{\sqrt{2}}\sinh(\sqrt{2}a_{0})+\left[\frac{a_{1}+b_{1}}{2}\cosh(\sqrt{2}a_{0})
+\frac{b_{1}-a_{1}}{2\sqrt{2}a_{0}}\sinh(\sqrt{2}a_{0})\right]n^{-\frac{1}{3}}\>+
\end{equation*}
\begin{equation*}
\left[\left(\frac{(a_{1}+b_{1})^{2}}{8a_{0}}-\frac{a_{1}^{2}}{2a_{0}}+\frac{a_{2}+b_{2}}{2}\right)\cosh(\sqrt{2}a_{0})
+\left(\frac{(a_{1}+b_{1})^{2}}{4\sqrt{2}}-\frac{(a_{1}+b_{1})^{2}}{8\sqrt{2}a_{0}^{2}}+\frac{a_{1}^{2}}{2\sqrt{2}a_{0}^{2}}\right.
\right.
\end{equation*}
\begin{equation}\label{matrix 31}
\left.\left.
+\>\>\frac{b_{2}-a_{2}}{2\sqrt{2}a_{0}}\right)\sinh(\sqrt{2}a_{0})\right]n^{-\frac{2}{3}}+\sinh(\sqrt{2}a_{0})O(n^{-1})
\>\>=\>\>\exp(M)_{31}\>\>=\>\>-\exp(M)_{23}.
\end{equation}
The last term of the exponential matrix \eqref{series solution} is
\begin{equation*}
1+\sum_{k\geq
1}\frac{2^{k}a^{k}b^{k}}{(2k)!}=\cosh(\sqrt{2}a_{0})+\frac{a_{1}+b_{1}}{\sqrt{2}}\sinh(\sqrt{2}a_{0})\,n^{-\frac{1}{3}}\>+\>
\left[\frac{(a_{1}+b_{1})^{2}}{4}\cosh(\sqrt{2}a_{0})\>+\right.
\end{equation*}
\begin{equation}\label{matrix 33}
\left.\left(\frac{a_{1}b_{1}}{\sqrt{2}a_{0}}+\frac{a_{2}+b_{2}}{\sqrt{2}}-\frac{(a_{1}+b_{1})^{2}}{4\sqrt{2}a_{0}}
\right) \sinh(\sqrt{2}a_{0})\right]n^{-\frac{2}{3}} +
\cosh(\sqrt{2}a_{0})\,O(n^{-1})\>\>=\>\>\exp(M)_{33}.
\end{equation}
We note that $\sqrt{2}a_{0}$ is exactly the quantity $\mu$ defined
in \cite{Trac2} as the new variable when solving for the limiting
system of equation as $n$ goes to infinity.

We then use equations \eqref{solution V_{n}}, \eqref{matrix 11},
\eqref{matrix 31} together with the numerical value of $c_{\varphi}$
for $n$ even given by \eqref{c varphi} to have the following
expansion of $V_{n,\varepsilon}$ and $q_{n,\varepsilon}$.
\begin{equation*}
V_{n,\varepsilon}(\tau(s))\>=\>
\frac{1}{2}(1-\exp(-\sqrt{2}a_{0}))+\left[\frac{a_{1}-b_{1}}{2\sqrt{2}}\frac{\sinh(\sqrt{2}a_{0})}{\sqrt{2}a_{0}}
-\frac{a_{1}+b_{1}}{2\sqrt{2}}\exp(-\sqrt{2}a_{0})\right]n^{-\frac{1}{3}}\>+
\end{equation*}
\begin{equation*}
\left(
\left[\frac{(a_{1}+b_{1})^{2}}{8}(1-\frac{1}{\sqrt{2}a_{0}})-\frac{\sqrt{2}(a_{2}+b_{2})}{4}+\frac{a_{1}^{2}}{2\sqrt{2}a_{0}}
\right]\cosh(\sqrt{2}a_{0})
-\left[\frac{(a_{1}+b_{1})^{2}}{8}(1+\frac{1}{\sqrt{2}a_{0}})\right.\right.
\end{equation*}
\begin{equation*}
\left. \left.
-\frac{(a_{1}+b_{1})^{2}}{16a_{0}^{2}}-\frac{a_{1}b_{1}}{2\sqrt{2}a_{0}}+\frac{a_{0}^{2}}{4a_{0}^{2}}
-\frac{\sqrt{2}(a_{2}-b_{2})}{4 \sqrt{2}a_{0}}
-\frac{\sqrt{2}(a_{2}+b_{2})}{4 }\right]\sinh(\sqrt{2}a_{0})
\right)n^{-\frac{2}{3}} +O(\frac{1}{n})
\end{equation*}
and
\begin{equation*}
q_{n,\varepsilon}(\tau(s))\>=\>
\frac{1}{\sqrt{2}}\exp(-\sqrt{2}a_{0})
+\left[-\frac{a_{1}+b_{1}}{2}\exp(-\sqrt{2}a_{0})
-\frac{a_{1}-b_{1}}{2}\frac{\sinh(\sqrt{2}a_{0})}{\sqrt{2}a_{0}}\right]n^{-\frac{1}{3}}
\end{equation*}
\begin{equation*}
\left(\left[\frac{\sqrt{2}(a_{1}+b_{1})^{2}}{8}(1-\frac{1}{\sqrt{2}a_{0}})+\frac{b_{1}^{2}}{2a_{0}}
-\frac{a_{2}+b_{2}}{2}\right]\cosh(\sqrt{2}a_{0}) +
\left[\frac{b_{2}-a_{2}}{2\sqrt{2}a_{0}} + \frac{a_{2}+b_{2}}{2}
\right.\right.
\end{equation*}
\begin{equation*}
\left.\left.-\frac{\sqrt{2}(a_{1}+b_{1})^{2}}{8}(1+\frac{1}{\sqrt{2}a_{0}})
+
\frac{a_{1}b_{1}}{2a_{0}}+\frac{\sqrt{2}(a_{1}+b_{1})^{2}}{16a_{0}^{2}}-
\frac{\sqrt{2}b_{1}^{2}}{4a_{0}^{2}}\right]\sinh(\sqrt{2}a_{0})
\right)n^{-\frac{2}{3}} +O(\frac{1}{n}).
\end{equation*}

\subsubsection{Second system of equations involving the calligraphic
variables for GOE$_{n}$}

The system involving the calligraphic variables is
\begin{equation}
Y(t)= \exp\left(-\int_{t}^{\iy}B(x)dx\right)\cdot Y(\iy).
\end{equation}
We set
\begin{displaymath}
\exp\left(-\int_{t}^{\iy}B(x)\,dx\right)= \exp\left\{\left(
               \begin{array}{ccc}
                 0 & 0 & -\int_{t}^{\iy}q_{n}(x)\,dx \\\\
                 0 & 0 & -\int_{t}^{\iy}p_{n}(x)\,dx \\ \\
                 -\int_{t}^{\iy}p_{n}(x) dx& -\int_{t}^{\iy}q_{n}(x)\,dx& 0 \\
               \end{array}
             \right)\right\}
\end{displaymath}
\begin{equation*}
=\exp(\mathcal{M}).
\end{equation*}
$\mathcal{M}$ is of the form
\begin{equation}
\mathcal{M}=\left(
               \begin{array}{ccc}
                 0 & 0 & -a \\
                 0 & 0 & -b \\
                 -b & -a & 0 \\
               \end{array}
             \right)
\end{equation}
with
\begin{displaymath}\label{call series solution}
\exp(\mathcal{M})=\left(
               \begin{array}{ccc}
                 1+\sum_{k\geq 1}\frac{2^{k-1}a^{k}b^{k}}{(2k)!} & \sum_{k\geq 1}\frac{2^{k-1}a^{k+1}b^{k-1}}{(2k)!} &
-\sum_{k\geq 0}\frac{2^{k}a^{k+1}b^{k}}{(2k+1)!}  \\ \\
                 \sum_{k\geq 1}\frac{2^{k-1}a^{k-1}b^{k+1}}{(2k)!}  & 1+\sum_{k\geq 1}\frac{2^{k-1}a^{k}b^{k}}{(2k)!}  &
-\sum_{k\geq 0}\frac{2^{k}a^{k}b^{k+1}}{(2k+1)!}  \\ \\
                 -\sum_{k\geq 0}\frac{2^{k}a^{k}b^{k+1}}{(2k+1)!}  & -\sum_{k\geq 0}\frac{2^{k}a^{k+1}b^{k}}{(2k+1)!}  &
1+\sum_{k\geq 1}\frac{2^{k}a^{k}b^{k}}{(2k)!}  \\
               \end{array}
             \right)
\end{displaymath}
\begin{equation}
=\left(\begin{array}{ccc} \exp(\mathcal{M})_{11} &\exp(\mathcal{M})_{12}&\exp(\mathcal{M})_{13}\\
\exp(\mathcal{M})_{21}&\exp(\mathcal{M})_{22}&\exp(\mathcal{M})_{23}\\
\exp(\mathcal{M})_{31}&\exp(\mathcal{M})_{32}&\exp(\mathcal{M})_{33}\\
\end{array}
\right)
\end{equation}
We have
\begin{equation}\label{call solution Q}
\mathcal{Q}_{n,1}(t)\>\>=\>\>
2c_{\varphi}\exp(\mathcal{M})_{11}\>\>+\>\>\exp(\mathcal{M})_{13},
\end{equation}
\begin{equation}\label{call solution P}
\mathcal{P}_{n,1}(t)\>\>=\>\>
2c_{\varphi}\exp(\mathcal{M})_{21}\>\>+\>\>\exp(\mathcal{M})_{23},
\end{equation}
and
\begin{equation}\label{call solution R}
\tilde{\mathcal{R}}_{n,1}(t)\>\>=\>\>
2c_{\varphi}\exp(\mathcal{M})_{31}\>\>+\>\>\exp(\mathcal{M})_{33}.
\end{equation}

We note that $\exp(\mathcal{M}_{21})=-\exp(M)_{21}, \>
\exp(\mathcal{M}_{23})=\exp(M)_{23},\,$
$\exp(\mathcal{M}_{31})=-\exp(M)_{31}\,$ and
$\,\exp(\mathcal{M}_{33})=\exp(M)_{33}$. The solutions \eqref{call
solution P} and \eqref{call solution R} follow directly from the
large $n$ expansion obtained in the last subsection. We therefore
have the following solutions for $\mathcal{P}_{n,1}$ and
$\tilde{\mathcal{R}}_{n,1}$

\begin{equation*}
\mathcal{P}_{n,1}(\tau(s))\>=\>
\frac{1}{\sqrt{2}}(\exp(-\sqrt{2}a_{0})-1) +\left[
\frac{a_{1}-b_{1}}{\sqrt{2}a_{0}}+\frac{a_{1}+b_{1}}{2}\exp(-\sqrt{2}a_{0})+\frac{b_{1}-a_{1}}{\sqrt{2}a_{0}}
\cosh(\sqrt{2}a_{0})\right.
\end{equation*}
\begin{equation*}
-\left.\frac{b_{1}-a_{1}}{\sqrt{2}a_{0}}\frac{\sinh(\sqrt{2}a_{0})}{2})\right]n^{-\frac{1}{3}}+\left(-\frac{\sqrt{2}a_{1}^{2}
}{2a_{0}^{2}}+\frac{\sqrt{2}a_{1}b_{1}}{2a_{0}^{2}}+\frac{a_{2}-b_{2}}{\sqrt{2}a_{0}}+\left[
\frac{(a_{1}+b_{1})^{2}}{4\sqrt{2}}+\frac{\sqrt{2}(a_{1}^{2}+a_{1}b_{1})}{2a_{0}^{2}}
\right..\right.
\end{equation*}
\begin{equation*}
\left.-\frac{a_{1}^{2}}{4a_{0}}-\frac{a_{2}-b_{2}}{\sqrt{2}a_{0}}+\frac{b_{1}^{2}}{4a_{0}}\right]\cosh(\sqrt{2}a_{0})+\left[
-\frac{(a_{1}+b_{1})^{2}}{4\sqrt{2}}+\frac{(a_{1}+b_{1})^{2}}{8\sqrt{2}a_{0}^{2}}-\frac{a_{1}^{2}}{2\sqrt{2}a_{0}^{2}}\right.
\end{equation*}
\begin{equation}
  \left.\left.
+\frac{a_{2}-b_{2}}{2\sqrt{2}a_{0}}\right]\sinh(\sqrt{2}a_{0})\right)n^{-\frac{2}{3}}
+O(\frac{1}{n})
\end{equation}
and
\begin{equation*}
\tilde{\mathcal{R}}_{n,1}(\tau(s))\>=\>\exp(-\sqrt{2}a_{0})+\left[\frac{a_{1}-b_{1}}{\sqrt{2}}\frac{\sinh(\sqrt{2}a_{0})}
{\sqrt{2}a_{0}}-\frac{a_{1}+b_{1}}{\sqrt{2}}\exp(-\sqrt{2}a_{0})\right]n^{-\frac{1}{3}}
\end{equation*}
\begin{equation*}
\left[\biggl(\frac{(a_{1}+b_{1})^{2}}{4}
-\frac{a_{2}+b_{2}}{\sqrt{2}}\biggr)\exp(-\sqrt{2}a_{0})-\frac{(a_{1}+b_{1})^{2}}{4\sqrt{2}a_{0}}\cosh(\sqrt{2}a_{0})-
\frac{(a_{1}-b_{1})^{2}}{4\sqrt{2}a_{0}}\sinh(\sqrt{2}a_{0})\right.
\end{equation*}
\begin{equation*}
\left.\frac{(a_{1}+b_{1})^{2}}{8a_{0}^{2}}\sinh(\sqrt{2}a_{0})+\frac{a_{2}-b_{2}}{2a_{0}}\sinh(\sqrt{2}a_{0})
+\frac{a_{1}^{2}}{\sqrt{2}a_{0}}\cosh(\sqrt{2}a_{0})-\frac{a_{1}^{2}}{2a_{0}^{2}}\sinh(\sqrt{2}a_{0})
\right]n^{-\frac{2}{3}}
\end{equation*}
\begin{equation}
+O(\frac{1}{n})
\end{equation}
\subsubsection{Calligraphic variables for GSE$_{n}$}
We note that the GSE$_{n}$ case is identical to the GOE$_{n}$ up to
a sign change and the parity of $n$ for the calligraphic variables.
The large $n$ expansion for $u_{n,\varepsilon}(t)$ and
$\tilde{v}_{n,\varepsilon}(t)$ follows from the matrix exponential
\eqref{series solution}. The boundary conditions need to be change
to $u_{n,\varepsilon}(\infty)=0$ and
$\tilde{v}_{n,\varepsilon}(\infty)=0$ and
$q_{n,\varepsilon}(\infty)=0$ as $n$ is odd.  Therefore in this case

\begin{equation}
u_{n,\varepsilon}(t)\>\>=\>\>
\exp(M)_{12}\>\>
\end{equation}
\begin{equation}
V_{n,\varepsilon}(t)\>\>=\>\>
\exp(M)_{22}\>\>
\end{equation}
and
\begin{equation}
q_{n,\varepsilon}(t)\>\>=\>\>
\exp(M)_{32}\>\>
\end{equation}
The large $n$ expansions of these quantities is given by
\eqref{matrix 12}, \eqref{matrix 11} and \eqref{matrix 13}
respectively.

The system of equations satisfied by the calligraphic variables is
\begin{equation}
\frac{d}{dt}\left(
  \begin{array}{c}
    \mathcal{Q}_{n,4}(t) \\
    \mathcal{P}_{n,4}(t) \\
    \tilde{\mathcal{R}}_{n,4}(t) \\
  \end{array}
\right) =\left(
               \begin{array}{ccc}
                 0 & 0 & -q_{n}(t) \\
                 0 & 0 & -p_{n}(t) \\
                 -p_{n}(t) & -q_{n}(t) & 0 \\
               \end{array}
             \right)\,\cdot\,
\left(
  \begin{array}{c}
    \mathcal{Q}_{n,4}(t) \\
    \mathcal{P}_{n,4}(t) \\
    \tilde{\mathcal{R}}_{n,4}(t) \\
  \end{array}
\right).
\end{equation}
where $ \tilde{\mathcal{R}}_{n,4}(t) =1+\mathcal{R}_{n,4}(t)$. The
boundary conditions in this case are
\begin{equation}
\left(
  \begin{array}{c}
    \mathcal{Q}_{n,4}(\infty) \\
    \mathcal{P}_{n,4}(\infty) \\
    \tilde{\mathcal{R}}_{n,4}(\infty) \\
  \end{array}
\right) =\left(
               \begin{array}{c}
                 -c_{\varphi} \\
                 -c_{\psi} \\
                1 \\
               \end{array}
             \right)\,=\,
\left(
  \begin{array}{c}
    0 \\
    -c_{\psi} \\
    1 \\
  \end{array}
\right)\quad \textrm{ as } n \textrm{ is odd}
\end{equation}
We can use the same technique as for the Orthogonal case to find a
large $n$ expansion of
$\mathcal{Q}_{n,4}(t),\>\>\mathcal{P}_{n,4}(t)$  and
    $\tilde{\mathcal{R}}_{n,4}(t) $.
\subsection{Large $n$ expansion of the probability distribution of
the largest eigenvalue}

We recall that the quantity of interest is \eqref{determinant}.
Under our change of variables it reads
\begin{equation}\label{result}
(1-\tilde{v}_{n,\varepsilon})(1-\frac{1}{2}\mathcal{R}_{n,1})-\frac{1}{2}(q_{n,\varepsilon}-c_{\varphi})\mathcal{P}_{n,1}=
\frac{1}{2}\left[V_{n,\varepsilon}\bigl(1+\tilde{\mathcal{R}}_{n,1}\bigr)-\mathcal{P}_{n,1}\bigl(q_{n,\varepsilon}-c_{\varphi}
\bigr)\right].
\end{equation}
Upon substitutions of the newly derived expressions in the right of
\eqref{result}, the right side of \eqref{result} takes the form
\begin{equation*}
\exp(-\sqrt{2}a_{0})
+\left[\frac{a_{1}-b_{1}}{2a_{0}}\bigl(1-\exp(-\sqrt{2}a_{0})\bigr)-\frac{a_{1}+b_{1}}{\sqrt{2}}\exp(-\sqrt{2}a_{0})\right]
n^{-\frac{1}{3}}
\end{equation*}
\begin{equation*}
\left\{\frac{a_{1}b_{1}-a_{1}^{2}}{2a_{0}^{2}}+\frac{a_{2}-b_{2}}{2a_{0}}+\left(\frac{3a_{1}^{2}-b_{1}^{2}}{8\sqrt{2}a_{0}}
-\frac{(a_{1}+b_{1})^{2}}{16\sqrt{2}a_{0}}-\frac{a_{2}+b_{2}}{4\sqrt{2}}\right)\exp(-2\sqrt{2}a_{0})\right.
\end{equation*}
\begin{equation*}
+\left(\frac{(a_{1}+b_{1})^{2}}{4}- \frac{3(a_{2}+b_{2})}{4\sqrt{2}}
+\frac{a_{1}^{2}}{2a_{0}^{2}}-
\frac{(a_{1}+b_{1})^{2}}{16\sqrt{2}a_{0}}+\frac{a_{1}^{2}}{2\sqrt{2}a_{0}}-\frac{a_{1}b_{1}}{4\sqrt{2}a_{0}}
-\frac{a_{2}-b_{2}}{2a_{0}}\right)\exp(-\sqrt{2}a_{0})
\end{equation*}
\begin{equation}\label{n expansion}
\left.
\frac{(a_{1}+b_{1})^{2}}{8a_{0}^{2}}\sinh(\sqrt{2}a_{0})-\frac{a_{1}b_{1}}{2a_{0}^{2}}\cosh(\sqrt{2}a_{0})\right\}
n^{-\frac{2}{3}} + O(\frac{1}{n})
\end{equation}
We follow Tracy and Widom \cite{Trac2} and denote by $\mu(s)$ the
quantity $\sqrt{2}a_{0}(s)$,
\begin{equation}
\mu:=\mu(s)\>=\>\sqrt{2}a_{0}(s)=\int_{s}^{\iy}q(x)dx,
\end{equation}
and we introduce the following notations,
\begin{equation}
\nu:=\nu(s)\>=\>\int_{s}^{\iy}p(x)dx,\quad \alpha:=\alpha(s)\>=\>
\int_{s}^{\iy}q(x)\,u(x)\,dx,
\end{equation}

\begin{equation*}
a_{1}(s)=\frac{1}{\sqrt{2}}\int_{s}^{\iy}\left(\frac{2c-1}{2}p(x)-cq(x)u(x)\right)dx,
\end{equation*}

\begin{equation*}
b_{1}(s)=\frac{1}{\sqrt{2}}\int_{s}^{\iy}\left(\frac{2c+1}{2}p(x)-cq(x)u(x)\right)dx,
\end{equation*}

\begin{equation*}
a_{2}(s)=\frac{1}{20\sqrt{2}}\int_{s}^{\iy}\left(\bigl(10c^{2}-10c+\frac{3}{2}\bigr)q_{1}(x)+p_{2}(x)+
\bigl(-30c^{2}+10c+\frac{3}{2}\bigr)q(x)v(x)+p_{1}v(x)\right.
\end{equation*}
\begin{equation*}
\left.+p(x)v_{1}(x)-q_{2}(x)u(x)-q_{1}(x)u_{1}(x)-q(x)u_{2}(x)+\bigl(\frac{3}{2}-10c^{2}\bigr)p(x)u(x)
+20c^{2}q(x)u^{2}(x)\right)\,dx,
\end{equation*}

\begin{equation*}
b_{2}(s)=\frac{1}{20\sqrt{2}}\int_{s}^{\iy}\left(\bigl(10c^{2}+10c+\frac{3}{2}\bigr)q_{1}(x)+p_{2}(x)+
\bigl(-30c^{2}-10c+\frac{3}{2}\bigr)q(x)v(x)+p_{1}v(x)\right.
\end{equation*}
\begin{equation*}
\left.+p(x)v_{1}(x)-q_{2}(x)u(x)-q_{1}(x)u_{1}(x)-q(x)u_{2}(x)+\bigl(\frac{3}{2}-10c^{2}\bigr)p(x)u(x)
+20c^{2}q(x)u^{2}(x)\right)\,dx.
\end{equation*}

We note that
\begin{equation}
a_{1}(s)-b_{1}(s)\>=\>-\frac{1}{\sqrt{2}}\int_{s}^{\iy}p(x)\,dx,
\end{equation}
\begin{equation}
a_{1}(s)+b_{1}(s)\>=\>\frac{1}{\sqrt{2}}\int_{s}^{\iy}\bigl(2cp(x)-2cq(x)u(x)\bigr)dx\>=\>-\sqrt{2}\,c\,q(s),
\end{equation}
and
\begin{equation}
a_{2}(s)-b_{2}(s)=\frac{c}{\sqrt{2}}\int_{s}^{\iy}\bigl(q(x)v(x)-q_{1}(x)\bigr)dx\>=\>\frac{c}{\sqrt{2}}p(s).
\end{equation}
We set
\begin{equation*}
\eta(s) \>=\> a_{2}(s)+
b_{2}(s)=\frac{1}{20\sqrt{2}}\int_{s}^{\iy}\biggl(\bigl(20c^{2}+3\bigr)q_{1}(x)+2p_{2}(x)+
\bigl(-60c^{2}+3\bigr)q(x)v(x)+2p_{1}v(x)
\end{equation*}
\begin{equation*}
+2p(x)v_{1}(x)-2q_{2}(x)u(x)-2q_{1}(x)u_{1}(x)-2q(x)u_{2}(x)+\bigl(3-20c^{2}\bigr)p(x)u(x)
\end{equation*}
\begin{equation}
+40c^{2}q(x)u^{2}(x)\biggr)\,dx.
\end{equation}
\begin{equation}
=\>
-\frac{20c^{2}q'(s)+3p(s)}{20\sqrt{2}}+\frac{1}{20\sqrt{2}}\int_{s}^{\infty}\bigl(6qv+3pu+2p_{2}+
2p_{1}v+2pv_{1}-2q_{2}u-2q_{1}u_{1}-2qu_{2}\bigr)(x)\,dx.
\end{equation}
This last equality comes from
\begin{equation*}(20c^{2}+3)q_{1}(x)+(3-60c^{2})q(x)v(x)+(3-20c^{2})p(x)u(x)+40c^{2}q(x)u^{2}(x)\>=\end{equation*}
\begin{equation}\label{eq8}
20c^{2}\bigl(q_{1}(x)-3q(x)v(x)-p(x)u(x)+2q(x)u^{2}(x)\bigr)+3\bigl(q_{1}(x)+q(x)v(x)+p(x)u(x)\bigr),
\end{equation}
the substitutions
\begin{equation}
q_{1}(x)=xq(x)-q(x)v(x)+p(x)u(x),\quad \textrm{and} \quad
q_{1}(x)=p'(x)+q(x)v(x) ,
\end{equation}
in \eqref{eq8}
\begin{equation*}
20c^{2}\bigl(xq(x)-4q(x)v(x)+2q(x)u^{2}(x)\bigr)+3\bigl(p'(x)+2q(x)v(x)+p(x)u(x)\bigr)
\end{equation*}
\begin{equation}\label{eq9}
=\>20c^{2}\bigl(xq(x)+2q(x)\bigl(-2v(x)+u^{2}(x)\bigr)\bigr)+3\bigl(p'(x)+2q(x)v(x)+p(x)u(x)\bigr),
\end{equation}
and the substitutions
\begin{equation}
q^{2}(x)=u^{2}(x)-2v(x),\quad \textrm{and}\quad
q''(x)=xq(x)+2q^{3}(x)
\end{equation}
in \eqref{eq9}
\begin{equation*}
20c^{2}\bigl(xq(x)+2q^{3}(x)\bigr)+3\bigl(p'(x)+2q(x)v(x)+p(x)u(x)\bigr)
\end{equation*}
\begin{equation}
=\>20c^{2}q''(x)+3p'(x) +6q(x)v(x)+3p(x)u(x).
\end{equation}

With these representations, \eqref{n expansion} is
\begin{equation*}
e^{-\mu}\>
+\>\left[c\,q(s)e^{-\mu}\>-\>\frac{1}{2\mu}\nu\bigl(1-e^{-\mu}\bigr)\right]
n^{-\frac{1}{3}}\>+\>\left\{\frac{2c-1}{4\mu^{2}}\nu^{2}-\frac{c}{2\mu^{2}}\nu\int_{s}^{\iy}q(x)u(x)dx
+\frac{c\,p(s)}{2\mu}\right.
\end{equation*}
\begin{equation*}
\left[\frac{c^{2}-2c+1/4}{8\mu}\nu^{2}-\frac{c(c-1)}{4\mu}\nu\,
\int_{s}^{\iy}q(x)u(x)dx +
\frac{c^{2}}{8\mu}\left(\int_{s}^{\iy}q(x)u(x)dx\right)^{2}
-\frac{c^{2}\,q^{2}(s)}{8\mu}\right.
\end{equation*}
\begin{equation*}
\left.-\frac{\eta}{4\sqrt{2}} \right]e^{-2\mu}\>+\>
\left[\frac{c^{2}\,q^{2}(s)}{2} -\frac{3 \eta}{4\sqrt{2}}
+\frac{2-\mu}{2\mu^{2}}\left(\frac{c^2-c+1/4}{2}\nu^{2}-\frac{2c^2-c}{2}\nu\int_{s}^{\iy}q(x)u(x)dx\right.\right.
\end{equation*}
\begin{equation*}
\left. +\frac{c^2}{2}\bigl(\int_{s}^{\iy}q(x)u(x)dx\bigr)^{2}\right)
-\frac{c^{2}\,q^{2}(s)}{8\mu}
-\frac{c^{2}-1/4}{8\mu}\nu^{2}+\frac{c^2}{4\mu}\nu\,
\int_{s}^{\iy}q(x)u(x)dx - \frac{c\,p(s)}{2\mu}
\end{equation*}
\begin{equation*}
\left.-\frac{c^{2}}{4\mu}\left(\int_{s}^{\iy}q(x)u(x)dx\right)^{2}\right]e^{-\mu}
+\frac{c^{2}\,q^{2}(s)}{2\mu^2}\sinh(\mu)
-\left(\frac{c^2-1/4}{2}\nu^2 -c^2 \nu \int_{s}^{\iy}q(x)u(x)dx
\right.
\end{equation*}
\begin{equation}\label{eq6}
\left. \left.c^2 \bigl(\int_{s}^{\iy}q(x)u(x)dx\bigr)^2
\right)\frac{\cosh(\mu)}{\mu^2} \right\}n^{-\frac{2}{3}}
\>+\>O(\frac{1}{n})
\end{equation}
as $n$ goes to infinity uniformly for $s$ bounded away from minus
infinity.\\

Finally we combine \eqref{eq6} with Theorem\ref{Fn2} to have the following version of Theorem\ref{GOE}.\\

If we set $t=\tau(s)$, then as $n\ra\iy$
\begin{equation*}
F_{n,1}^{2}(t)\>=\>F_{2}(s)\cdot\left\{e^{-\mu}\>+\>
\left[c\bigl(q(s)+u(s)\bigr) e^{-\mu } -\frac{\nu }{2 \mu
}(1-e^{-\mu }) \right]n^{-\frac{1}{3}} \>\> + \right.
\end{equation*}
\begin{equation*}
\left[-\frac{1}{20} E_{c,2}(s)\,e^{-\mu } -\frac{c\,\alpha(s)}{2 \mu
^2}+\frac{c\, p(s)}{2 \mu }+\frac{(2c-1)\, \nu ^2}{4 \mu ^2}
+c\,u(s) \left(c\,q(s)\,e^{-\mu } -\frac{\nu }{2 \mu }(1-e^{-\mu }
)\right) \right.
\end{equation*}
\begin{equation*}
+\>\>e^{-2 \mu } \left(-\frac{\eta }{4
\sqrt{2}}+\frac{c^2\,\alpha^2(s) }{8 \mu }-\frac{c^2\, q^2(s)}{8 \mu
}-\frac{ (c^2-c) \, \nu \,\alpha(s)}{4 \mu
}+\frac{\left(\frac{1}{4}-2 c+c^2\right)\, \nu ^2}{8 \mu
}\right)\>\>+
\end{equation*}
\begin{equation*}
e^{-\mu } \left(\frac{c^2 \,q^2(s)}{2}-\frac{3\, \eta }{4
\sqrt{2}}-\frac{c^2 \,\alpha^2(s) }{4 \mu }-\frac{c\, p(s)}{2 \mu
}-\frac{c^2\, q^2(s)}{8 \mu }+\frac{ c^2 \, \nu \,\alpha(s)}{4 \mu
}-\frac{\left(-\frac{1}{4}+c^2\right) \nu^2}{8 \mu }+\right.
\end{equation*}
\begin{equation*}
\left. \frac{2-\mu }{2 \mu ^2} \left(\frac{ c^2\,\alpha^2(s)}{2}
-\frac{  (2c^{2}-c)\,\nu \,\alpha(s) }{2} +\frac{
\left(\frac{1}{4}-c+c^2\right)\, \nu^2}{2}\right)\right) \>\>-
\end{equation*}
\begin{equation}\label{eq7}
\left.\left.\left( c^2\,\alpha^2(s)- c^2 \,\nu\, \alpha(s) +\frac{
\left(-\frac{1}{4}+c^2\right)\, \nu ^2}{2}\right)\frac{ \cosh(\mu)
}{\mu ^2}+\frac{c^2\, q^2(s) }{8 \mu ^2}\sinh(\mu)\right]
n^{-\frac{2}{3}}\>\right\}\>+\>O(n^{-1})
\end{equation}
uniformly in $s$.\\

To simplify the $n^{-\frac{2}{3}}$ term in \eqref{eq7}, we use the
representation $p(s)=q'(s)+q(s)u(s)$ which says in this setting that
$ \nu(s)=\alpha(s)-q(s)$. The result of this substitution is
Theorem\ref{GOE}.

\section{Conclusion}
We note that unlike $F_{n,2}(t)$ for the GUE$_{n}$, the GOE$_{n}$
large $n$ expansion of the probability distribution of the largest
eigenvalue $F_{n,1}(t)$ has a non vanishing $n^{-\frac{1}{3}}$
correction term. Thus the convergence to the limiting Tracy-Widom
distribution $F_{1}(t)$ is slower. Numerical applications of
$F_{n,1}(t)$ follows easily from $q(s)$ this is one consequence of
our representation of $F_{n,1}(t)$ in Theorem\ref{GOE}. All the
terms on the right side of \eqref{Fn1} can be expressed in terms of
$q(s)$ and $q'(s)$. \\
The GSE$_{n}$ largest eigenvalue distribution is derived in a
similar way (the only major difference being that $n$ needs to be
odd in this case.)

\clearpage \vspace{3ex} \noindent\textbf{\large Acknowledgements: }
The author would like to thank Professor Craig Tracy for the
discussions that initiated this work and for the invaluable
guidance, and the Department of Mathematical Sciences at the
University of Alabama in Huntsville.


\begin{thebibliography}{10}

\bibitem{Deif1}
J.~Baik, P.~A.~Deift and  K.~Johansson.
\newblock { On the distribution of the length of the
longest increasing subsequence in a random permutation}
\newblock{J. Amer. Math. Soc., 12 (1999), 1119–1178.}

\bibitem{Choup1}
L.~N.~Choup.
\newblock{Edgeworth Expansion of the Largest Eigenvalue Distribution Function
 of GUE and LUE }
\newblock{\em IMRN }Volume 2006, ID 61049, Pages 1-33.
\bibitem{Choup2}
L.~N.~Choup.
\newblock { Edgeworth Expansion of the Largest Eigenvalue Distribution Function of GUE
Revisited}
\newblock ArXiv:0711.4206v1

\bibitem{Deif2}
P.~Deift.
\newblock{Orthogonal Polynomials and Random Matrices: A Riemann-Hilbert Approach}.
\newblock {\em American Mathematical Society}. Courant Lecture
Notes 3, 2000.

\bibitem{Deif3}
P.~Deift,
\newblock{Universality for mathematical and physical systems.}
\newblock{\em International Congress of Mathematicians,} Vol.1, 125-152, Eur.Math.Soc., Z$\ddot{u}$rich, 2007.

\bibitem{Dien1}
M.~Dieng and C.~A.~Tracy.
\newblock{Application of random matrix theory to multivariate statistics}.
\newblock{preprint, Arxiv:math.PR/0603543}.


\bibitem{Fell2}
W.~Feller.
\newblock {\em An Introduction to Probability Theory and Its
Applications} ,Vol.II.
\newblock Second edition, John Wiley, 1971.

\bibitem{Forr1}
T.~M.~Garoni, P.~J.~Forrester and N.~E.~Frankel.
\newblock{Asymptotic corrections to the eigenvalue density of the
GUE and LUE.}
\newblock arXiv:math-ph/0504053 v1

\bibitem{Gohb1}
I.~Gohberg, S.~Goldberg, and M.~A. Kaashoek.
\newblock {\em {Classes of Linear Operators, Vol. I}}, volume~49 of {\em
  Operator Theory: Advances and Applications}.
\newblock Birkh{\"a}user, 1990.

\bibitem{Gohb2}
I.~Gohberg, S.~Goldberg, and M.~A. Kaashoek.
\newblock {\em {Classes of Linear Operators, Vol. II}}, volume~63 of {\em
  Operator Theory: Advances and Applications}.
\newblock Birkh{\"a}user, 1993.

\bibitem{Gohb3}
I.~C. Gohberg, M.~G. Kre$\breve{i}$n.
\newblock {\em {Introduction to the Theory of Linear Nonselfadjoint Operators}}, volume~18 of {\em
  Translations of Mathematical Monographs}.
\newblock American Mathematical Society, 1969.



\bibitem{Hoch1}
H.~Hochstadt.
\newblock {\em {The Functions of Mathematical Physics}}, volume~23 of {\em
  Pure and Applied Mathematics: A series of texts and Monographs }.
\newblock Wiley-Interscience, 1971.

\bibitem{Johansson}
K.~Johansson.
\newblock{ Toeplitz determinants, random growth and determinantal
processes.}
\newblock{\em Proceedings of the ICM, Beijing 2002,} vol. 3,
53--62, math.PR/0304368.

\bibitem{John1}
I.~M.~Johnstone,
\newblock {On the distribution of the largest eigenvalue in principal component
  analysis},
\newblock {\em Ann. Stats.}, 29(2):295--327, 2001.


\bibitem{Lax1}
P.~D.~Lax.
\newblock{ \em Functional Analysis}
\newblock{ Wiley-Interscience}, 2002.


\bibitem{Meht1}
M.~L.~Mehta.
\newblock {\em Random Matrices, Revised and Enlarged Second Edition}.
\newblock Academic Press, 1991.

\bibitem{Olve1}
F.~W.~J.~Olver.
\newblock{Asymptotics and Special Functions}
\newblock Academic Press, New York, 1974.

\bibitem{Plan1}
M.~Plancherel and W. Rotach.
\newblock{ Sur les valeurs asymptotiques des polynomes d'Hermite}
\newblock Comm. Math. Helv. 1 (1929)227-254.

\bibitem{Sosh1}
A.~Soshnikov.
\newblock {Universality at the Edge of the Spectrum in Wigner Ranom Matrices}.
\newblock {\em J. Stat. Phys.}, 108(5--6):1033--1056, 2002.

\bibitem{Szeg1}
G.~Szeg\"o.
\newblock{Orthogonal Polynomials.}
\newblock American Mathematical Society Colloquium Publications
Volume 23

\bibitem{Taylor}
M.~E.~Taylor
\newblock{Partial Differential Equations}.
\newblock{Springer-Verlag, New York, 1996}

\bibitem{Trac3}
C.~A.~Tracy and H.~Widom.
\newblock {Level--spacing distributions and the Airy kernel}.
\newblock {\em Commun. Math. Physics}, 159:151--174, 1994.

\bibitem{Trac7}
C.~A.~Tracy and H.~Widom.
\newblock {Fredholm determinants, differential equations and matrix models}.
\newblock {\em Commun. Math. Physics}, 163:33--72, 1994.

\bibitem{Trac2}
C.~A.~Tracy and H.~Widom.
\newblock {On orthogonal and symplectic matrix ensembles}.
\newblock {\em Commun. Math. Physics}, 177:727--754, 1996.

\bibitem{Trac1}
C.~A.~Tracy and H.~Widom.
\newblock {Correlation functions, cluster functions, and spacing distributions
  for random matrices}.
\newblock {\em J. Stat. Phys.}, 92(5--6):809--835, 1998.

\bibitem{Trac4}
C.~A. Tracy and H.~Widom.
\newblock {Airy kernel and Painlev\'e II}.
\newblock In {\em Isomonodromic deformations and applications in physics},
  volume~31 of {\em {CRM Proceedings \& Lecture Notes}}, pages 85--98. Amer.
  Math. Soc., Providence, RI, 2002.

\bibitem{Trac8}
C.~A.~Tracy and H.~Widom.
\newblock {Distribution functions for largest eigenvalues
and their applications}.
\newblock In {\em Proceedings of the International Congress
of Mathematicians, Beijing 2002}, Vol.~I, ed. LI Tatsien, Higher
Education Press, Beijing, pgs.~587--596, 2002.

\bibitem{Trac5}
C.~A.~Tracy and H.~Widom.
\newblock {Matrix kernels for the Gaussian orthogonal and symplectic
  ensembles}.
\newblock {\em Ann. Inst. Fourier, Grenoble}, 55, 2197--2207, 2005.

\bibitem{Whit1}
E.~T.~Whittaker and G.~N.~Watson.
\newblock{\em A Course of Modern Analysis} Fourth Edition
\newblock{Cambridge University Press}, 2004.


\end{thebibliography}
\end{document}